\documentclass{article}
\usepackage{amssymb}

\input{tcilatex}
\begin{document}

\begin{center}
\bigskip \textbf{Morgan type uncertainty principle} \textbf{and unique
continuation properties for abstract Schr\"{o}dinger equations}

\textbf{Veli Shakhmurov}

Department of Mechanical Engineering, Okan University, Akfirat, Tuzla 34959
Istanbul, Turkey,

E-mail: veli.sahmurov@okan.edu.tr

A\textbf{bstract}
\end{center}

In this paper, Morgan type uncertainty principle and unique continuation
properties of abstract Schr\"{o}dinger equations with time dependent
potentials in vector-valued $L^{2}$ classes are obtained. The equation
involves a possible linear operators considered in the Hilbert space $H.$
So, by choosing the corresponding spaces $H$ and operators we derived unique
continuation properties for numerous classes of Schr\"{o}dinger type
equations and its systems which occur in a wide variety of physical systems.


\begin{center}
\bigskip\ \ \textbf{AMS 2010: 35Q41, 35K15, 47B25, 47Dxx, 46E40 }

\textbf{1. Introduction, definitions}
\end{center}

\bigskip In this paper, the unique continuation properties of the abstract
Schr\"{o}dinger equations

\begin{equation}
i\partial _{t}u+\Delta u+Au+V\left( x,t\right) u=0,\text{ }x\in R^{n},\text{ 
}t\in \left[ 0,1\right] ,  \tag{1.1}
\end{equation}%
are studied, where $A$ is a linear operator$,$ $V\left( x,t\right) $ is a
given potential operator function in a Hilbert space $H$, subscript $t$
indicates the partial derivative with respect to $t$, $\Delta $ denotes the
Laplace operator in $R^{n}$ and $u=$ $u(x,t)$ is the $H$-valued unknown
function. The goal is to obtain sufficient conditions on the opertator $A$,
the potential $V$ and the behavior of the solution $u$ at two different
times, $t_{0}=0$ and $t_{1}=1$ which guarantee that $u\equiv 0$ in $%
R^{n}\times \lbrack 0,1].$

This linear result was then applied to show that two regular solutions $%
u_{1} $ and $u_{2}$ of non-linear Schr\"{o}dinger equations 
\begin{equation}
i\partial _{t}u+\Delta u+Au=F\left( u,\bar{u}\right) ,\text{ }x\in R^{n},%
\text{ }t\in \left[ 0,1\right] ,  \tag{1.2}
\end{equation}%
and for general non-linearities $F$, must agree in $R^{n}\times \lbrack 0,1]$%
, when $u_{1}-u_{2}$ and its gradient decay faster than any quadratic
exponential at times $0$ and $1$.

Unique continuation properties for Schr\"{o}dinger equations studied e.g in $%
\left[ \text{6-8}\right] ,$ $\left[ \text{21-23}\right] $ and the referances
therein.\ In contrast to the mentioned above results we will study the
unique continuation properties of abstract Schr\"{o}dinger equations with
operator potentials. Abstract differential equations studied e.g. in $\left[
1\right] $, $\left[ 5\right] $, $\left[ 10\text{, }11\right] ,\left[ 13-19%
\right] $, $\left[ 25\text{, }26\right] .$ Since the Hilbert space $H$ is
arbitrary and $A$ is a possible linear operator, by choosing $H$ and $A$ we
can obtain numerous classes of Schr\"{o}dinger type equations and its
systems which occur in a wide variety of physical systems. If we\ choose the
abstract space $H$ a concrete Hilbert space, for example $H=L^{2}\left(
G\right) $, $A=L,$ where $G$ is a domin in $R^{m}$ with sufficientli smooth
boundary and $L$ is elliptic operator, then we obtain the unique
continuation properties of followinng Schr\"{o}dinger equation%
\begin{equation}
\partial _{t}u=i\left[ \Delta u+Lu+V\left( x,t\right) u\right] ,\text{ }x\in
R^{n},\text{ }y\in G,\text{ }t\in \left[ 0,1\right] ,  \tag{1.3}
\end{equation}%
where $L$ is an elliptic operator with respect to variable $y\in G\subset
R^{m}$ and $u=u\left( x,y,t\right) $.

\ Moreover, let $H=L^{2}\left( 0,1\right) $ and $A$ to be differential
operator with generalized Wentzell-Robin boundary condition defined by 
\begin{equation}
D\left( A\right) =\left\{ u\in W^{2,2}\left( 0,1\right) ,\text{ }%
B_{j}u=Au\left( j\right) =0\text{, }j=0,1\right\} ,\text{ }  \tag{1.4}
\end{equation}%
\[
\text{ }Au=au^{\left( 2\right) }+bu^{\left( 1\right) },
\]%
where $a$ is positive and $b$ is a real-valued functions on $\left(
0,1\right) $. Then, we get the unique continuation properties of the
Wentzell-Robin type boundary value problem (BVP) for the following Schr\"{o}%
dinger type equation 
\begin{equation}
i\partial _{t}u+\Delta u+a\frac{\partial ^{2}u}{\partial y^{2}}+b\frac{%
\partial u}{\partial y}=F\left( u,\bar{u}\right) ,\text{ }  \tag{1.5}
\end{equation}%
\[
u=u\left( x,y,t\right) \text{, }x\in R^{n},\text{ }y\in G,\text{ }t\in \left[
0,1\right] ,
\]%
\ \ \ 

\begin{equation}
a\left( j\right) u_{yy}\left( x,j,t\right) +b\left( j\right) u_{y}\left(
x,j,t\right) =0\text{, }j=0,1\text{, for a.e. }x\in R^{n},\text{ }t\in
\left( 0,1\right) .  \tag{1.6}
\end{equation}

Note that, the regularity properties of Wentzell-Robin type BVP for elliptic
equations were studied e.g. in $\left[ \text{11, 12 }\right] $ and the
references therein. Moreover, if put $H=l_{2}$ and choose $A$ and $V\left(
x,t\right) $ as a infinite matrices $\left[ a_{mj}\right] $, $b_{mj}\left(
x,t\right) $ $m,$ $j=1,2,...,\infty $ respectively, then we obtain the
unique continuation properties of the following system of Schr\"{o}dinger
equation 
\begin{equation}
\partial _{t}u_{m}=i\left[ \Delta u_{m}+\sum\limits_{j=1}^{N}\left(
a_{mj}+b_{mj}\left( x,t\right) \right) u_{j}\right] ,\text{ }x\in R^{n},%
\text{ }t\in \left( 0,1\right) .  \tag{1.7}
\end{equation}

Let $E$ be a Banach space and $\gamma =\gamma \left( x\right) $ be a
positive measurable function on a domain $\Omega \subset R^{n}.$Here, $%
L_{p,\gamma }\left( \Omega ;E\right) $ denotes the space of strongly
measurable $E$-valued functions that are defined on $\Omega $ with the norm

\[
\left\Vert f\right\Vert _{_{p,\gamma }}=\left\Vert f\right\Vert
_{L_{p,\gamma }\left( \Omega ;E\right) }=\left( \int \left\Vert f\left(
x\right) \right\Vert _{E}^{p}\gamma \left( x\right) dx\right) ^{\frac{1}{p}},%
\text{ }1\leq p<\infty , 
\]

\[
\left\Vert f\right\Vert _{L_{\infty ,\gamma }\left( \Omega ;E\right)
}=ess\sup\limits_{x\in \Omega }\left\Vert f\left( x\right) \right\Vert
_{E}\gamma \left( x\right) ,\text{ }p=\infty . 
\]

For $\gamma \left( x\right) \equiv 1$ the space $L_{p,\gamma }\left( \Omega
;E\right) $ will be denoted by $L_{p}=L_{p}\left( \Omega ;E\right) $ for $%
p\in \left[ 1,\infty \right] .$

Let $\mathbf{p=}\left( p_{1},p_{2}\right) $ and $\Omega =\Omega _{1}\times
\Omega _{2}$, where $\Omega _{k}\in R^{n_{k}}$. $L_{x}^{p_{1}}L_{t}^{p_{2}}%
\left( \Omega ;E\right) $ will denote the space of all $E$-valed $\mathbf{p}$%
-summable\ functions with mixed norm, i.e., the space of all measurable
functions $f$ defined on $G$ equiped with norm%
\[
\left\Vert f\right\Vert _{L_{x}^{p_{1}}L_{y}^{p_{2}}\left( \Omega ;E\right)
}=\left( \left( \dint\limits_{\Omega _{1}}\left\Vert f\left( x,t\right)
\right\Vert _{E}^{p_{2}}dt\right) ^{\frac{p_{1}}{p_{2}}}dx\right) ^{\frac{1}{%
p_{1}}}. 
\]

For $p=2$ and $H$ Hilbert space we get Hilbert space of $H$-valued functions
with inner product of two elements $f$, $g\in L^{2}\left( \Omega ;H\right) $:%
\[
\left( f,g\right) _{L^{2}\left( \Omega ;H\right) }=\int\limits_{\Omega
}\left( f\left( x\right) ,g\left( x\right) \right) _{H}dx. 
\]

Let $C\left( \Omega ;E\right) $ denote the space of $E-$valued, bounded
uniformly continious functions on $\Omega $ with norm 
\[
\left\Vert u\right\Vert _{C\left( \Omega ;E\right) }=\sup\limits_{x\in
\Omega }\left\Vert u\left( x\right) \right\Vert _{E}. 
\]

$C^{m}\left( \Omega ;E\right) $\ will denote the space of $E$-valued bounded
uniformly strongly continuous and $m$-times continuously differentiable
functions on $\Omega $ with norm 
\[
\left\Vert u\right\Vert _{C^{m}\left( \Omega ;E\right) }=\max\limits_{0\leq
\left\vert \alpha \right\vert \leq m}\sup\limits_{x\in \Omega }\left\Vert
D^{\alpha }u\left( x\right) \right\Vert _{E}. 
\]

$C_{0}^{\infty }\left( \Omega ;E\right) -$denotes the space of $E$-valued
infinity many differentiable finite functions.

Let

\[
O_{R}=\left\{ x\in R^{n},\text{ }\left\vert x\right\vert <R\right\} \text{, }%
R>0. 
\]
Let $\mathbb{N}$ denote the set of all natural numbers, $\mathbb{C}$ denote
the set of all complex numbers.

Let $\Omega $ be a domain in $R^{n}$ and $m$ be a positive integer$.$\ $%
W^{m,p}\left( \Omega ;E\right) $ denotes the space of all functions $u\in
L^{p}\left( \Omega ;E\right) $ that have the generalized derivatives $\frac{%
\partial ^{m}u}{\partial x_{k}^{m}}\in L^{p}\left( \Omega ;E\right) ,$ $%
1\leq p\leq \infty $ with the norm 
\[
\ \left\Vert u\right\Vert _{W^{m,p}\left( \Omega ;E\right) }=\left\Vert
u\right\Vert _{L^{p}\left( \Omega ;E\right)
}+\sum\limits_{k=1}^{n}\left\Vert \frac{\partial ^{m}u}{\partial x_{k}^{m}}%
\right\Vert _{L^{p}\left( \Omega ;E\right) }<\infty . 
\]

Let $E_{0}$ and $E$ be two Banach spaces and $E_{0}$ is continuously and
densely embedded into $E$. Here,\ $W^{m,p}\left( \Omega ;E_{0},E\right) $
denote the space $W^{m,p}\left( \Omega ;E\right) \cap $ $L^{p}\left( \Omega
;E\right) $ equipped with norm 
\[
\ \left\Vert u\right\Vert _{W^{m,p}\left( \Omega ;E_{0,}E\right)
}=\left\Vert u\right\Vert _{L^{p}\left( \Omega ;E_{0}\right)
}+\sum\limits_{k=1}^{n}\left\Vert \frac{\partial ^{m}u}{\partial x_{k}^{m}}%
\right\Vert _{L^{p}\left( \Omega ;E\right) }<\infty . 
\]

Let $E_{1}$ and $E_{2}$ be two Banach spaces. $L\left( E_{1},E_{2}\right) $
will denote the space of all bounded linear operators from $E_{1}$ to $%
E_{2}. $ For $E_{1}=E_{2}=E$ it will be denoted by $L\left( E\right) .$

A linear operator\ $A$ is said to be positive in a Banach\ space $E$ with
bound $M>0$ if $D\left( A\right) $ is dense on $E$ and $\left\Vert \left(
A+sI\right) ^{-1}\right\Vert _{L\left( E\right) }\leq M\left( 1+\left\vert
s\right\vert \right) ^{-1}$ for any $s\in $ $\left( -\infty ,0\right) ,$
where $I$ is the identity operator in $E.$

Let $\left[ A,B\right] $ be a commutator operator, i.e. 
\[
\left[ A,B\right] =AB-BA 
\]%
for linear operators $A$ and $B.$

Sometimes we use one and the same symbol $C$ without distinction in order to
denote positive constants which may differ from each other even in a single
context. When we want to specify the dependence of such a constant on a
parameter, say $\alpha $, we write $C_{\alpha }$.

\begin{center}
\textbf{2}. \textbf{Main results for} \textbf{absract Scr\"{o}dinger equation%
}
\end{center}

First of all, we generalize the result G. W. Morgan (see e.g $\left[ 7\right]
$) about Morgan type uncertainty principle for Fourier transform. Let 
\[
X=L^{2}\left( R^{n};H\right) \text{, }Y^{k}=W^{2,k}\left( R^{n};H\right) 
\text{, }k\in \mathbb{N}. 
\]

\textbf{Lemma 2.1. }Let \textbf{\ }$f\left( x\right) \in L^{1}\left(
R^{n};H\right) \cap X$ and%
\[
\dint\limits_{R^{n}}\dint\limits_{R^{n}}\left\Vert f\left( x\right)
\right\Vert _{H}\left\Vert \hat{f}\left( \xi \right) \right\Vert
_{H}e^{\left\vert x.\zeta \right\vert }\text{ }dxd\xi <\infty \text{.} 
\]

Then $f\left( x\right) \equiv 0.$

In particular, using Young's inequality this implies:

\textbf{Result 2.1.} Let%
\[
f\left( x\right) \in L^{1}\left( R^{n};H\right) \cap X,\text{ }p\in \left(
1,2\right) ,\frac{1}{p}+\frac{1}{q}=1,\alpha ,\beta >0 
\]%
and 
\[
\dint\limits_{R^{n}}\left\Vert f\left( x\right) \right\Vert _{H}e^{\frac{%
\alpha ^{p}\left\vert x\right\vert ^{p}}{p}}dx+\dint\limits_{R^{n}}\left%
\Vert \hat{f}\left( \xi \right) \right\Vert _{H}e^{\frac{\beta
^{q}\left\vert \xi \right\vert ^{q}}{q}}d\xi <\infty \text{, }\alpha \beta
>1. 
\]

Then $f\left( x\right) \equiv 0.$

The Morgan type uncertainty principle, in terms of the solution of the free
Schr\"{o}dinger equation will be as:

Let%
\[
u_{0}\left( .\right) \in L^{1}\left( R^{n};H\right) \cap X 
\]%
and for some $t\neq 0$ 
\[
\dint\limits_{R^{n}}\left\Vert u_{0}\left( x\right) \right\Vert _{H}e^{\frac{%
\alpha ^{p}\left\vert x\right\vert ^{p}}{p}}dx+\dint\limits_{R^{n}}\left%
\Vert e^{it\left( \Delta +A\right) }u_{0}\left( x\right) \right\Vert _{H}e^{%
\frac{\beta ^{q}\left\vert \xi \right\vert }{q\left( 2t\right) ^{q}}%
^{q}}dx<\infty \text{, }\alpha \beta >1. 
\]

Then $u_{0}\left( x\right) \equiv 0.$

\textbf{Condition 1. }Assume $A$ is a pozitive operator in Hilbert space $H$
and $iA$ generates a semigrop $U\left( t\right) =e^{iAt}$. \ Suppose%
\begin{equation}
\left\Vert V\right\Vert _{L^{\infty }\left( R^{n}\times \left( 0,1\right)
;L\left( H\right) \right) }\leq C  \tag{2.1}
\end{equation}

and 
\begin{equation}
\text{ }\lim\limits_{R\rightarrow \infty }\left\Vert V\right\Vert
_{L_{t}^{1}L_{x}^{\infty }\left( L\left( H\right) \right) }=0,  \tag{2.2}
\end{equation}

where 
\[
\text{ }L_{t}^{1}L_{x}^{\infty }\left( L\left( H\right) \right) =L^{1}\left(
0,1;L^{\infty }\left( R^{n}/O_{R}\right) ;L\left( H\right) \right) . 
\]

Here, 
\[
\sigma \left( t\right) =\frac{1}{\alpha \left( 1-t\right) +\beta t}. 
\]%
In $\left[ \text{17, Theorem 1}\right] $ we proved the following result:

\textbf{Theorem A}$_{1}$\textbf{. }Assume the Condition 1 holds or $V\left(
x,t\right) =V_{1}\left( x\right) +V_{2}\left( x,t\right) $, where $%
V_{1}\left( x\right) \in L\left( H\right) $ for $x\in R^{n}$ and%
\[
\sup\limits_{t\in \left[ 0,1\right] }\left\Vert e^{\left\vert x\right\vert
^{2}\sigma ^{2}\left( t\right) }V_{2}\left( .,t\right) \right\Vert
_{B}<\infty . 
\]

Suppose $\alpha $, $\beta >0$ and $\alpha \beta <4$ such that any solution $%
u\in C\left( \left[ 0,1\right] ;X\right) $ of $\left( 1.1\right) $ satisfay 
\[
\left\Vert e^{\frac{\left\vert x\right\vert ^{2}}{\beta ^{2}}}u\left(
.,0\right) \right\Vert _{X}<\infty ,\left\Vert e^{\frac{\left\vert
x\right\vert ^{2}}{\alpha ^{2}}}u\left( .,1\right) \right\Vert _{X}<\infty . 
\]

Then $u\left( x,t\right) \equiv 0.$

Our main result in this paper is the following:

\textbf{Theorem 1. }Assume the Condition 1 holds and there exist constants $%
a_{0},$ $a_{1},$ $a_{2}>0$ such that for any $k\in \mathbb{Z}^{+}$ a
solution $u\in C\left( \left[ 0,1\right] ;X\right) $ of $\left( 1.1\right) $
satisfy 
\begin{equation}
\dint\limits_{R^{n}}\left\Vert u\left( x,0\right) \right\Vert
_{H}^{2}e^{2a_{0}\left\vert x\right\vert ^{p}}dx<\infty ,\text{ for }p\in
\left( 1,2\right) ,  \tag{2.3}
\end{equation}%
\begin{equation}
\dint\limits_{R^{n}}\left\Vert u\left( x,1\right) \right\Vert
_{H}^{2}e^{2k\left\vert x\right\vert ^{p}}dx<a_{2}e^{2a_{1}k^{\frac{q}{q-p}%
}},\text{ }\frac{1}{p}+\frac{1}{q}=1.  \tag{2.4}
\end{equation}%
Moreover, there exists $M_{p}>0$ such that

\begin{equation}
a_{0}a_{1}^{p-2}>M_{p}.  \tag{2.5}
\end{equation}

Then $u\left( x,t\right) \equiv 0.$

\bigskip \textbf{Corollary 1. }Assume the Condition 1 holds and 
\[
\lim\limits_{\left\vert R\right\vert \rightarrow \infty
}\dint\limits_{0}^{1}\sup\limits_{\left\vert x\right\vert >R}\left\Vert
V\left( x,t\right) \right\Vert _{L\left( H\right) }dt=0. 
\]%
There exist positive constants $\alpha $, $\beta $ such that a solution $%
u\in C\left( \left[ 0,1\right] ;X\right) $ of $\left( 1.1\right) $ satisfy 
\begin{equation}
\dint\limits_{R^{n}}\left\Vert u\left( x,0\right) \right\Vert
_{H}^{2}e^{2\left\vert \alpha x\right\vert
^{p}/p}dx+\dint\limits_{R^{n}}\left\Vert u\left( x,1\right) \right\Vert
_{H}^{2}e^{\frac{2\left\vert \beta x\right\vert ^{q}}{q}}dx<\infty ,\text{ }
\tag{2.6}
\end{equation}%
with 
\[
p\in \left( 1,2\right) \text{, }\frac{1}{p}+\frac{1}{q}=1 
\]
and there exists $N_{p}>0$ such that 
\begin{equation}
\alpha \beta >N_{p}.\text{ }  \tag{2.7}
\end{equation}%
Then $u\left( x,t\right) \equiv 0.$

As a direct consequence of Corollary 1 we have the following result
regarding the uniqueness of solutions for nonlinear equation $\left(
1.2\right) $

\bigskip \textbf{Theorem 2. }Assume the Condition 1 holds and $u_{1},$ $%
u_{2}\in C\left( \left[ 0,1\right] ;Y^{2,k}\right) $ strong solutions of $%
(1.2)$ with $k\in \mathbb{Z}^{+},$ $k>\frac{n}{2}.$ Suppose $F:H\times
H\rightarrow H,$ $F\in C^{k}$, $F\left( 0\right) =\partial _{u}F\left(
0\right) =\partial _{\bar{u}}F\left( 0\right) =0$ and there exist positive
constants $\alpha $, $\beta $ such that%
\begin{equation}
e^{\frac{\left\vert \alpha x\right\vert ^{p}}{p}}\left( u_{1}\left(
.,0\right) -u_{2}\left( .,0\right) \right) \in X,\text{ }e^{^{\frac{%
\left\vert \beta x\right\vert ^{q}}{q}}}\left( u_{1}\left( .,0\right)
-u_{2}\left( .,0\right) \right) \in X,  \tag{2.8}
\end{equation}%
with 
\[
p\in \left( 1,2\right) \text{, }\frac{1}{p}+\frac{1}{q}=1 
\]
and there exists $N_{p}>0$ such that 
\begin{equation}
\alpha \beta >N_{p}.\text{ }  \tag{2.9}
\end{equation}

Then $u_{1}\left( x,t\right) \equiv u_{2}\left( x,t\right) .$

\textbf{Corollary 2.} Assume the Condition 1 holds and there exist positive
constants $\alpha $ and $\beta $ such that a solution $u\in C\left( \left[
0,1\right] ;X\right) $ of $\left( 1.1\right) $ satisfy 
\begin{equation}
\dint\limits_{R^{n}}\left\Vert u\left( x,0\right) \right\Vert _{H}^{2}e^{%
\frac{2\left\vert \alpha x_{j}\right\vert ^{p}}{p}}dx+\dint\limits_{R^{n}}%
\left\Vert u\left( x,1\right) \right\Vert _{H}^{2}e^{\frac{2\left\vert \beta
x_{j}\right\vert ^{q}}{q}}dx<\infty ,\text{ }  \tag{2.10}
\end{equation}%
for $j=1,2,...,n$ and $p\in \left( 1,2\right) $, $\frac{1}{p}+\frac{1}{q}=1.$
Moreover, there exists $N_{p}>0$ such that 
\begin{equation}
\alpha \beta >N_{p}.\text{ }  \tag{2.11}
\end{equation}%
Then $u\left( x,t\right) \equiv 0.$

\textbf{Remark 2.1. }The Theorem 3 still holds, with different constant $%
N_{p}>0$, if one replaces the hypothesis $(2.8)$ by

\[
e^{\left\vert \alpha x_{j}\right\vert ^{p}/p}\left( u_{1}\left( .,0\right)
-u_{2}\left( .,0\right) \right) \in X,\text{ }e^{\left\vert \beta
x_{j}\right\vert ^{q}/q}\left( u_{1}\left( .,0\right) -u_{2}\left(
.,0\right) \right) \in X, 
\]

\bigskip for $j=1,2,...,n.$

Next, we shall extend the method used in the proof Theorem 3 to study the
blow up phenomenon of solutions of nonlinear Schr\"{o}dinger equations

\begin{equation}
i\partial _{t}u+\Delta u+Au+F\left( u,\bar{u}\right) u=0,\text{ }x\in R^{n},%
\text{ }t\in \left[ 0,1\right] ,  \tag{2.12}
\end{equation}%
\bigskip where $A$ is a linear operatorin a Hilbert space $H.$

Let $u\left( x,t\right) $ be a solution of the equation $\left( 2.12\right)
. $ Then it can be shown that the function%
\begin{equation}
\upsilon \left( x,t\right) =U\left( x,1-t\right) u\left( \frac{x}{1-t},\frac{%
t}{1-t}\right) \text{, }  \tag{2.13}
\end{equation}%
is a solution of the focussing $L^{2}$-critical solution of abstract
Schrodinger equation%
\begin{equation}
i\partial _{t}u+\Delta u+Au+\left\Vert u\right\Vert ^{\frac{4}{n}}u=0,\text{ 
}x\in R^{n},\text{ }t\in \left[ 0,1\right]  \tag{2.14}
\end{equation}%
which blows up at time $t=1,$ where $U\left( x,t\right) $ is a fundamental
solution of the Schr\"{o}dinger equation%
\[
i\partial _{t}u+\Delta u+Au=0,\text{ }x\in R^{n},\text{ }t\in \left[ 0,1%
\right] , 
\]
i.e. 
\[
U\left( x,t\right) =\frac{1}{t^{n/2}}\exp \left\{ i\left( A+\left\vert
x\right\vert ^{2}\right) /4t\right\} . 
\]
Assume that the functon $\upsilon $ def\i ned by $\left( 2\text{.}13\right) $
satisfies the following estimate 
\begin{equation}
\left\Vert \upsilon \left( x,t\right) \right\Vert _{H}\leq \frac{1}{\left(
1-t\right) ^{n/2}}Q\left( \frac{\left\vert x\right\vert }{1-t}\right) ,\text{
}t\in \left( -1,1\right)  \tag{2.15}
\end{equation}%
where 
\begin{equation}
Q\left( x\right) =b_{1}^{-\frac{n}{2}}e^{-b_{2}\left\vert x\right\vert ^{p}}%
\text{, }b_{1}\text{, }b_{2}>0\text{, }p>1  \tag{2.16}
\end{equation}

The following result will be get:

\textbf{Theorem 3. }Assume the Condition 1 holds and there exist positive
constants $b_{0}$ and $\theta $ such that a solution $u\in C\left( \left[
-1,1\right] ;X\right) $ of $\left( 2.12\right) $ satisfied 
\begin{equation}
\left\Vert F\left( u,\bar{u}\right) \right\Vert _{H}\leq b_{0}\left\Vert
u\right\Vert _{H}^{\theta }\text{ for }\left\Vert u\right\Vert _{H}>1. 
\tag{2.17}
\end{equation}

Suppose 
\[
\left\Vert u\left( .,t\right) \right\Vert _{X}=\left\Vert u\left( .,0\right)
\right\Vert _{X}=\left\Vert u_{0}\right\Vert _{X}=a,\text{ }t\in \left(
-1,1\right) 
\]%
and that $\left( 2.15\right) $ holds with $Q\left( .\right) $ satisfies $%
\left( 2.16\right) .$\ If $p>p\left( \theta \right) =\frac{2\left( \theta
n-2\right) }{\left( \theta n-1\right) },$ then $a\equiv 0.$

\begin{center}
\textbf{3. Some properties of solutions of abstract Schredinger equations}
\end{center}

Let

\begin{center}
\[
\sigma \left( t\right) =\left[ \alpha \left( 1-t\right) +\beta t\right]
^{-1},\text{ }\eta \left( x,t\right) =\left( \alpha -\beta \right) )|x|^{2}%
\left[ 4i(\alpha (1-t)+\beta t)\right] ^{-1},\text{ } 
\]
\end{center}

\[
\nu \left( s\right) =\left[ \gamma \alpha \beta \sigma ^{2}\left( s\right) +%
\frac{\left( \alpha -\beta \right) a}{4\left( a^{2}+b^{2}\right) }\sigma
\left( s\right) \right] ,\text{ }\phi \left( x,t\right) =\text{ }\frac{%
\gamma a\left\vert x\right\vert ^{2}}{a+4\gamma \left( a^{2}+b^{2}\right) t}%
. 
\]

We recall the following lemma (see $\left[ \text{17, Lemma 3.1}\right] $).

\bigskip Let 
\[
\Phi \left( A,V\right) \upsilon =a\func{Re}\left( \left( A+V\right) \upsilon
,\upsilon \right) _{H}-b\func{Im}\left( \left( A+V\right) \upsilon ,\upsilon
\right) _{H}, 
\]%
\[
\text{ for }\upsilon =\upsilon \left( x,t\right) \in H\left( A\right) . 
\]

\textbf{Lemma A}$_{1}$\textbf{. }Assume $a>0,$ $b\in \mathbb{R}$, $A$ is a
symmetric operator in $H.$ Moreover, there is a constant $C_{0}>0$ so that 
\[
\left\vert \Phi \left( A,V\right) \upsilon \left( x,t\right) \right\vert
\leq C_{0}\mu \left( x,t\right) \left\Vert \upsilon \left( x,t\right)
\right\Vert _{H}^{2}, 
\]%
for $x\in R^{n},$ $t\in \left[ 0,T\right] ,$ $\gamma \geq 0$, $T\in \left[
0,1\right] $ and $\upsilon \in H\left( A\right) $, where $\mu $ is a
positive function in $L^{1}\left( 0,T;L^{\infty }\left( R^{n}\right) \right) 
$.

Then the solution $u$ of $\left( 3.0\right) $ belonging to $L^{\infty
}\left( 0,1;X\right) \cap L^{2}\left( 0,1;Y^{1}\right) $ satisfies the
following estimate%
\[
e^{M_{T}}\left\Vert e^{\phi \left( .,T\right) }u\left( .,T\right)
\right\Vert _{X}\leq M_{T}\left\Vert e^{\gamma \left\vert x\right\vert
^{2}}u\left( .,0\right) \right\Vert _{X}+\sqrt{a^{2}+b^{2}}\left\Vert
e^{\phi \left( t\right) }F\right\Vert _{L^{1}\left( 0,T;X\right) }, 
\]%
where 
\[
\text{ }M_{T}=\left\Vert \mu \right\Vert _{L^{1}\left( 0,T:L^{\infty }\left(
R^{n}\right) \right) }. 
\]

Let $u=u\left( x,s\right) $ be a solution of the equation 
\[
\partial _{s}u=i\left[ \Delta u+Au+V\left( y,s\right) u+F\left( y,s\right) %
\right] ,\text{ }y\in R^{n},\text{ }s\in \left[ 0,1\right] . 
\]%
and $a+ib\neq 0$, $\gamma \in \mathbb{R}$, $\alpha $, $\beta \in \mathbb{R}%
_{+}$. Set 
\begin{equation}
\tilde{u}\left( x,t\right) =\left( \sqrt{\alpha \beta }\sigma \left(
t\right) \right) ^{\frac{n}{2}}u\left( \sqrt{\alpha \beta }x\sigma \left(
t\right) ,\beta t\sigma \left( t\right) \right) e^{\eta }.  \tag{3.1}
\end{equation}

\bigskip Then, $\tilde{u}\left( x,t\right) $\ verifies the equation 
\begin{equation}
\partial _{t}\tilde{u}=i\left[ \Delta \tilde{u}+A\tilde{u}+\tilde{V}\left(
x,t\right) \tilde{u}+\tilde{F}\left( x,t\right) \right] ,\text{ }x\in R^{n},%
\text{ }t\in \left[ 0,1\right]  \tag{3.2}
\end{equation}%
with 
\begin{equation}
\tilde{V}\left( x,t\right) =\alpha \beta \sigma ^{2}\left( t\right) V\left( 
\sqrt{\alpha \beta }x\sigma \left( t\right) ,\beta t\sigma \left( t\right)
\right) ,  \tag{3.3}
\end{equation}

\begin{equation}
\text{ }\tilde{F}\left( x,t\right) =\left( \sqrt{\alpha \beta }\sigma \left(
t\right) \right) ^{\frac{n}{2}+2}\left( \sqrt{\alpha \beta }x\sigma \left(
t\right) ,\beta t\sigma \left( t\right) \right) .  \tag{3.4}
\end{equation}%
Moreover, 
\begin{equation}
\left\Vert e^{\gamma \left\vert x\right\vert ^{2}}\tilde{F}\left( .,t\right)
\right\Vert _{X}=\alpha \beta \sigma ^{2}\left( t\right) e^{\nu \left\vert
y\right\vert ^{2}}\left\Vert F\left( s\right) \right\Vert _{X}\text{ and }%
\left\Vert e^{\gamma \left\vert x\right\vert ^{2}}\tilde{u}\left( .,t\right)
\right\Vert _{X}=e^{\nu \left\vert y\right\vert ^{2}}\left\Vert u\left(
s\right) \right\Vert _{X}  \tag{3.5}
\end{equation}%
when $s=\beta t\sigma \left( t\right) $.

\textbf{Remark 3.1. }Let $\beta =\beta \left( k\right) .$ By assumption we
have 
\[
\left\Vert e^{a_{0}\left\vert x\right\vert ^{p}}u\left( x,0\right)
\right\Vert _{X}=A_{0}, 
\]%
\begin{equation}
\left\Vert e^{k\left\vert x\right\vert ^{p}}u\left( x,0\right) \right\Vert
_{X}=A_{k}\leq a_{2}e^{2a_{1}k^{\frac{q}{q-p}}}=a_{2}e^{2a_{1}k^{\frac{1}{2-p%
}}}.  \tag{3.6}
\end{equation}

\bigskip Thus, for $\gamma =\gamma (k)$ $\in \lbrack 0,\infty )$ to be
chosen later, one has%
\begin{equation}
\left\Vert e^{\gamma \left\vert x\right\vert ^{p}}\tilde{u}_{k}\left(
x,0\right) \right\Vert _{X}=\left\Vert e^{\gamma \left( \frac{\alpha }{\beta 
}\right) ^{p/2}\left\vert x\right\vert ^{p}}u_{k}\left( x,0\right)
\right\Vert _{X}=B_{0},  \tag{3.7}
\end{equation}

\[
\left\Vert e^{\gamma \left\vert x\right\vert ^{p}}\tilde{u}_{k}\left(
x,1\right) \right\Vert _{X}=\left\Vert e^{\gamma \left( \frac{\beta }{\alpha 
}\right) ^{p/2}\left\vert x\right\vert ^{p}}u_{k}\left( x,1\right)
\right\Vert _{X}=A_{k}. 
\]

Let we take 
\[
\gamma \left( \frac{\alpha }{\beta }\right) ^{p/2}=a_{0}\text{ and }\gamma
\left( \frac{\beta }{\alpha }\right) ^{p/2}=k, 
\]

i.e. 
\begin{equation}
\gamma =\left( ka_{0}\right) ^{\frac{1}{2}}\text{, }\beta =k^{\frac{1}{p}},%
\text{ }\alpha =a_{0}^{\frac{1}{p}}.  \tag{3.8}
\end{equation}

Let 
\[
M=\dint\limits_{0}^{1}\left\Vert V\left( .,t\right) \right\Vert _{L^{\infty
}\left( R^{n};H\right) }dt=\dint\limits_{0}^{1}\left\Vert V\left( .,s\right)
\right\Vert _{L^{\infty }\left( R^{n};H\right) }ds. 
\]

From $\left( 3.2\right) $, using energy estimates it follows 
\begin{equation}
e^{-M}\left\Vert u\left( .,0\right) \right\Vert _{X}\leq \left\Vert u\left(
.,t\right) \right\Vert _{X}=\left\Vert \tilde{u}\left( .,s\right)
\right\Vert _{X}\leq e^{M}\left\Vert u\left( .,0\right) \right\Vert _{X}%
\text{, }t,s\in \left[ 0,1\right] ,  \tag{3.9}
\end{equation}%
where 
\[
s=\beta t\sigma \left( t\right) . 
\]

Consider the following problem%
\begin{equation}
i\partial _{t}u+\Delta u+Au=V\left( x,t\right) u+F\left( x,t\right) ,\text{ }%
x\in R^{n},\text{ }t\in \left[ 0,1\right] ,  \tag{3.10}
\end{equation}

\[
u\left( x,0\right) =u_{0}\left( x\right) ,\text{ } 
\]%
where $A$ is a linear operator$,$ $V\left( x,t\right) $ is a given potential
operator function in a Hilbert space $H$ and $F$ is a $H$-valued function.

Let as define operator valued integral operators in $L^{p}\left( \Omega
;E\right) $. Let $k$: $R^{n}\backslash \left\{ 0\right\} \rightarrow L\left(
E\right) .$ We say $k\left( x\right) $ is a $L\left( E\right) $-valued
Calderon-Zygmund kernel ($C-Z$ kernel) if $k\in C^{\infty }\left(
R^{n}\backslash \left\{ 0\right\} ,L\left( E\right) \right) ,$ $k$ is
homogenous of degree $-n,$ $\dint\limits_{B}k\left( x\right) d\sigma =0,$
where%
\[
B=\left\{ x\in R^{n}:\left\vert x\right\vert =1\right\} . 
\]
For $f\in L^{p}\left( \Omega ;E\right) ,$ $p\in \left( 1,\infty \right) ,$ $%
a\in L^{\infty }\left( R^{n}\right) ,$ and $x\in \Omega $ we set the
Calderon-Zygmund operator 
\[
K_{\varepsilon }f=\dint\limits_{\left\vert x-y\right\vert >\varepsilon ,y\in
\Omega }k\left( x,y\right) f\left( y\right) dy,\text{ }Kf=\lim\limits_{%
\varepsilon \rightarrow 0}K_{\varepsilon }f 
\]%
and commutator operator%
\[
\left[ K;a\right] f=a\left( x\right) Kf\left( x\right) -K\left( af\right)
\left( x\right) =\lim\limits_{\varepsilon \rightarrow
0}\dint\limits_{\left\vert x-y\right\vert >\varepsilon ,y\in \Omega }k\left(
x,y\right) \left[ a\left( x\right) -a\left( y\right) \right] f\left(
y\right) dy. 
\]

By using Calder\'{o}n's first commutator estimates $\left[ 4\right] ,$
convolution operators on abstract functions $\left[ 2\right] $ and abstract
commutaror theorem in $\left[ 20\right] $ we obtain the following result:

\textbf{Theorem A}$_{2}.$ Assume $k(.)$ is $L\left( E\right) $-valued $C-Z$
kernel that have locally integrable first-order derivatives in $\left\vert
x\right\vert >0$, and 
\[
\left\Vert k\left( x,y\right) -k\left( x^{\prime },y\right) \right\Vert
_{L\left( E\right) }\leq M\left\vert x-x^{\prime }\right\vert \left\vert
x-y\right\vert ^{-\left( n+1\right) }\text{ for }\left\vert x-y\right\vert
>2\left\vert x-x^{\prime }\right\vert .
\]%
Let $a(.)$ have first-order derivatives in $L^{r}\left( R^{n}\right) $,$%
1<r\leq \infty $. Then for $p$, $q\in \left( 1,\infty \right) ,$ $%
q^{-1}=p^{-1}+r^{-1}$ the following estimates hold%
\[
\left\Vert \left[ K;a\right] \partial _{x_{j}}f\right\Vert _{L^{q}\left(
R^{n};E\right) }\leq C\left\Vert f\right\Vert _{L^{p}\left( R^{n};E\right) },
\]%
\[
\left\Vert \partial _{x_{j}}\left[ K;a\right] f\right\Vert _{L^{q}\left(
R^{n};E\right) }\leq C\left\Vert f\right\Vert _{L^{p}\left( R^{n};E\right) },
\]%
for $f\in C_{0}^{\infty }\left( R^{n};E\right) ,$\ where the constant $C>0$
is independent of $f.$

Let \ 
\[
X_{\gamma }=L_{\gamma }^{2}\left( R^{n};H\right) . 
\]%
By following $\left[ \text{4, Lemma 2.1}\right] $ let us show the following

\textbf{Lemma 3.1.} Assume the Condition 1 holds and there exists $%
\varepsilon _{0}>0$ such that%
\[
\left\Vert V\right\Vert _{L_{t}^{1}L_{x}^{\infty }\left( R^{n}\times \left[
0,1\right] ;L\left( H\right) \right) }<\varepsilon _{0}. 
\]%
Moreover, suppose $u\in C\left( \left[ 0,1\right] ;X\right) $ is a stronge
solution of $\left( 3.10\right) $ with 
\[
u_{0},\text{ }u_{1}=u\left( x,1\right) \in X_{\gamma },\text{ }F\in
L^{1}\left( 0,1;X_{\gamma }\right) 
\]%
for $\gamma \left( x\right) =e^{2\lambda .x}$ $\ $and for some $\lambda \in
R^{n}.$ Then there exists a positive constant $M_{0}=M_{0}\left(
n,A,H\right) $ independent of $\lambda $ such that%
\begin{equation}
\sup\limits_{t\in \left[ 0,1\right] }\left\Vert u\left( .,t\right)
\right\Vert _{X_{\gamma }}\leq M_{0}\left[ \left\Vert u_{0}\right\Vert
_{X_{\gamma }}+\left\Vert u_{1}\right\Vert _{X_{\gamma
}}+\dint\limits_{0}^{1}\left\Vert F\left( .,t\right) \right\Vert _{X_{\gamma
}}dt\right] .  \tag{3.11}
\end{equation}

\textbf{Proof. }First, we consider the case, when $\gamma \left( x\right)
=\beta \left( x\right) =e^{2\beta x_{1}}.$ Without loss of generality we
shall assume $\beta >0.$ Let $\varphi _{n}\in C^{\infty }\left( \mathbb{R}%
\right) $ such that $\varphi _{n}\left( \tau \right) =1$, $\tau \leq n$ and $%
\varphi _{n}\left( \tau \right) =0$ for $\tau \geq 10n$ with $0\leq \varphi
_{n}\leq 1,$ $\left\vert \varphi _{n}^{\left( j\right) }\left( \tau \right)
\right\vert \leq C_{j}\tau n^{-j}.$ Let 
\[
\theta _{n}(\tau )=\beta \dint\limits_{0}^{\tau }\varphi _{n}^{2}\left(
s\right) ds 
\]%
so that $\theta _{n}\in C^{\infty }\left( \mathbb{R}\right) $ nondecreasing
with $\theta _{n}(\tau )=\beta \tau $ for $\tau <n,$ $\theta _{n}(\tau
)=C_{n}\beta $ for $\tau >10n$ and 
\begin{equation}
\theta _{n}^{\prime }(\tau )=\beta \varphi _{n}^{2}\left( \tau \right) \leq
\beta ,\text{ }\theta _{n}^{\left( j\right) }(\tau )=\beta C_{j}n^{1-j},%
\text{ }j=1,2,....  \tag{3.12}
\end{equation}

\bigskip Let $\phi _{n}\left( \tau \right) =\exp \left( 2\theta _{n}(\tau
)\right) $ so that $\phi _{n}\left( \tau \right) \leq \exp \left( 2\beta
\tau \right) $ and $\phi _{n}\left( \tau \right) \rightarrow \exp \left(
2\beta \tau \right) $ for $n\rightarrow \infty .$ Let $u\left( x,t\right) $
be a solution of the equation $\left( 3.10\right) $, then one gets the
equation $\upsilon _{n}\left( x,t\right) =\phi _{n}\left( x_{1}\right)
u\left( x,t\right) $ satisfies the following 
\begin{equation}
i\partial _{t}\upsilon _{n}+\Delta \upsilon _{n}+A\upsilon _{n}=V_{n}\left(
x,t\right) \upsilon _{n}+\phi _{n}\left( x_{1}\right) F\left( x,t\right) , 
\tag{3.13}
\end{equation}%
where 
\[
V_{n}\left( x,t\right) \upsilon _{n}=V\left( x,t\right) \upsilon _{n}+4\beta
\varphi _{n}\left( x_{1}\right) \partial _{x_{1}}\upsilon _{n}+\left[ 4\beta
\varphi _{n}\left( x_{1}\right) \varphi _{n}^{\prime }\left( x_{1}\right)
-4\beta ^{2}\upsilon _{n}^{4}\right] \upsilon _{n}. 
\]

Now, we consider a new function

\begin{equation}
w_{n}\left( x,t\right) =e^{\mu }\upsilon _{n}\left( x,t\right) ,\text{ }\mu
=-i4\beta ^{2}\varphi _{n}^{4}\left( x_{1}\right) t.  \tag{3.14}
\end{equation}

\bigskip Then from $\left( 3.13\right) $ we get%
\begin{equation}
i\partial _{t}w_{n}+\Delta w_{n}+Aw_{n}=\tilde{V}_{n}\left( x,t\right) w_{n}+%
\tilde{F}_{n}\left( x,t\right) ,  \tag{3.13}
\end{equation}%
where%
\[
\tilde{V}_{n}\left( x,t\right) w_{n}=V\left( x,t\right) w_{n}+h\left(
x_{1},t\right) +a^{2}\left( x_{1}\right) \partial _{x_{1}}w_{n}+itb\left(
x_{1}\right) \partial _{x_{1}}w_{n} 
\]%
\[
\tilde{F}_{n}\left( x,t\right) =e^{\mu }\phi _{n}\left( x_{1}\right) F\left(
x,t\right) , 
\]%
when 
\[
h\left( x_{1},t\right) =\left( i16\beta ^{2}\varphi _{n}^{3}\varphi
_{n}^{\prime }t\right) ^{2}+i48\beta ^{2}\varphi _{n}^{2}\left( \varphi
_{n}^{\prime }\right) ^{2}t+i16\beta ^{2}\varphi _{n}^{3}\varphi
_{n}^{\left( 2\right) }t+ 
\]

\[
4\beta \varphi _{n}\varphi _{n}^{\prime }+i64\beta ^{2}\varphi
_{n}^{3}\varphi _{n}^{\prime }t,\text{ }a^{2}=4\beta \varphi _{n}^{2}\left(
x_{1}\right) \text{, }b=-32\beta ^{2}\varphi _{n}^{3}\varphi _{n}^{\prime }. 
\]

It is clear to see that 
\begin{equation}
\left\Vert \partial _{x_{1}}^{j}h\left( x_{1},t\right) \right\Vert
_{L^{\infty }\left( \mathbb{R}\times \left[ 0,1\right] \right) }\leq
C_{j}n^{-\left( j+1\right) },\text{ }j=1,2,...,  \tag{3.14}
\end{equation}

\begin{equation}
a^{2}\left( x_{1}\right) \geq 0,\text{ }\left\Vert \partial
_{x_{1}}^{j}a\left( x_{1}\right) \right\Vert _{L^{\infty }\left( \mathbb{R}%
\right) }\leq C_{j}n^{-j},\text{ }j=1,2,...,  \tag{3.15}
\end{equation}

\[
\left\Vert \partial _{x_{1}}^{j}b\left( x_{1}\right) \right\Vert _{L^{\infty
}\left( \mathbb{R}\right) }\leq C_{j}n^{-j},\text{ }j=1,2,.... 
\]

Consider the smooth function $\eta \in C^{\infty }\left( R^{n}\right) $ such
that $\eta \left( x\right) =1$ for $\left\vert x\right\vert \leq \frac{1}{2}$
and $\eta \left( x\right) =0$ for $\left\vert x\right\vert \geq 1,$ with $%
0\leq \eta \left( x\right) \leq 1.$ Moreover, let 
\[
\chi _{\pm }\left( \xi _{1}\right) =\left\{ 
\begin{array}{c}
1\text{, }\xi _{1}>0(\xi _{1}<0) \\ 
0,\text{ }\xi _{1}<0\left( \xi _{1}>0\right)%
\end{array}%
\right. . 
\]

Define the multipliers operators%
\[
P_{\varepsilon }f=F^{-1}\left( \eta \left( \varepsilon \xi \right) \hat{f}%
\left( \xi \right) \right) ,\text{ }P_{\pm }f=F^{-1}\left( \chi _{\pm
}\left( \varepsilon \xi \right) \hat{f}\left( \xi \right) \right) \text{ for 
}0<\varepsilon \leq 1. 
\]

\bigskip Then by applying the equation $\left( 3.10\right) $ for example to $%
P_{\varepsilon }P_{+}w_{n}$ we get%
\begin{equation}
i\partial _{t}P_{\varepsilon }P_{+}w_{n}+\Delta P_{\varepsilon
}P_{+}w_{n}+AP_{\varepsilon }P_{+}w_{n}=P_{\varepsilon }P_{+}\left(
Vw_{n}\right) +P_{\varepsilon }P_{+}\left( hw_{n}\right) +\text{ } 
\tag{3.16}
\end{equation}

\[
P_{\varepsilon }P_{+}\left( a^{2}\left( x_{1}\right) \partial
_{x_{1}}w_{n}\right) +P_{\varepsilon }P_{+}\left( ib\left( x_{1}\right)
\partial _{x_{1}}w_{n}\right) +P_{\varepsilon }P_{+}\left( \tilde{F}%
_{n}\right) . 
\]

From $\left( 3.16\right) $ we obtain 
\[
i\partial _{t}\left( P_{\varepsilon }P_{+}w_{n},\upsilon \right) +\left(
\Delta P_{\varepsilon }P_{+}w_{n},\upsilon \right) +\left( AP_{\varepsilon
}P_{+}w_{n},\upsilon \right) =\text{ } 
\]

\[
\left( P_{\varepsilon }P_{+}\left( Vw_{n}\right) ,\upsilon \right) +\left(
P_{\varepsilon }P_{+}\left( hw_{n}\right) ,\upsilon \right) +\left(
P_{\varepsilon }P_{+}\left( a^{2}\left( x_{1}\right) \partial
_{x_{1}}w_{n}\right) ,\upsilon \right) + 
\]%
\begin{equation}
\left( P_{\varepsilon }P_{+}\left( ib\left( x_{1}\right) \partial
_{x_{1}}w_{n}\right) ,\upsilon \right) +\left( P_{\varepsilon }P_{+}\left( 
\tilde{F}_{n}\right) ,\upsilon \right) ,  \tag{3.17}
\end{equation}%
for all $\upsilon \in C_{0}^{\infty }\left( R^{n};H\right) $, where $\left(
u,\upsilon \right) $ denotes scalar product of $u\left( x,t\right) $ and $%
\upsilon \left( x\right) $ in $H$ for all $x\in R^{n}$ and $t\in \left[ 0,1%
\right] .$ Taking the complex conjugate of the above, we get the equation%
\[
\bar{K}\left( i\partial _{t}\left( P_{\varepsilon }P_{+}w_{n},\upsilon
\right) \right) +\bar{K}\left( \Delta P_{\varepsilon }P_{+}w_{n},\upsilon
\right) +\bar{K}\left( AP_{\varepsilon }P_{+}w_{n},\upsilon \right) =\text{ }
\]

\begin{equation}
\bar{K}\left( P_{\varepsilon }P_{+}\left( Vw_{n}\right) ,\upsilon \right) +%
\bar{K}\left( P_{\varepsilon }P_{+}\left( hw_{n}\right) ,\upsilon \right) +%
\bar{K}\left( P_{\varepsilon }P_{+}\left( a^{2}\left( x_{1}\right) \partial
_{x_{1}}w_{n}\right) ,\upsilon \right) +  \tag{3.18}
\end{equation}%
\[
\bar{K}\left( P_{\varepsilon }P_{+}\left( ib\left( x_{1}\right) \partial
_{x_{1}}w_{n}\right) ,\upsilon \right) +\bar{K}\left( P_{\varepsilon
}P_{+}\left( \tilde{F}_{n}\right) ,\upsilon \right) , 
\]%
here $\bar{K}\left( u\right) $ denotes the complex conjugate of $u.$

Multiplying $\left( 3.17\right) $ and $\left( 3.18\right) $ by $\bar{K}%
\left( \left( P_{\varepsilon }P_{+}w_{n},\upsilon \right) \right) $ and $%
-\left( P_{\varepsilon }P_{+}w_{n},\upsilon \right) $ respectively, and
adding the result, we obtain

\[
i\partial _{t}\left\vert \left( P_{\varepsilon }P_{+}w_{n},\upsilon \right)
\right\vert ^{2}+\left( \Delta P_{\varepsilon }P_{+}w_{n},\upsilon \right) 
\bar{K}\left( \left( P_{\varepsilon }P_{+}w_{n},\upsilon \right) \right) -%
\text{ } 
\]

\begin{equation}
\bar{K}\left( \Delta P_{\varepsilon }P_{+}w_{n},\upsilon \right) \left(
\left( P_{\varepsilon }P_{+}w_{n},\upsilon \right) \right) +\left(
AP_{\varepsilon }P_{+}w_{n},\upsilon \right) \bar{K}\left( \left(
P_{\varepsilon }P_{+}w_{n},\upsilon \right) \right) -  \tag{3.19}
\end{equation}

\[
\left( \left( P_{\varepsilon }P_{+}w_{n},\upsilon \right) \right) \bar{K}%
\left( AP_{\varepsilon }P_{+}w_{n},\upsilon \right) =\left( P_{\varepsilon
}P_{+}w_{n},\upsilon \right) \bar{K}\left( P_{\varepsilon }P_{+}\left(
Vw_{n}\right) ,\upsilon \right) - 
\]%
\[
\bar{K}\left( \left( P_{\varepsilon }P_{+}w_{n},\upsilon \right) \right)
\left( P_{\varepsilon }P_{+}\left( Vw_{n}\right) ,\upsilon \right) +\left(
\left( P_{\varepsilon }P_{+}w_{n},\upsilon \right) \right) \bar{K}\left(
P_{\varepsilon }P_{+}\left( \tilde{F}_{n}\right) ,\upsilon \right) - 
\]

\[
\bar{K}\left( P_{\varepsilon }P_{+}\left( \tilde{F}_{n}\right) ,\upsilon
\right) \left( P_{\varepsilon }P_{+}w_{n},\upsilon \right) +\left( \left(
P_{\varepsilon }P_{+}hw_{n}\right) ,\upsilon \right) \bar{K}\left( \left(
P_{\varepsilon }P_{+}w_{n},\upsilon \right) \right) - 
\]%
\[
\bar{K}\left( \left( P_{\varepsilon }P_{+}hw_{n}\right) ,\upsilon \right)
\left( P_{\varepsilon }P_{+}w_{n},\upsilon \right) +\left( \left(
P_{\varepsilon }P_{+}a^{2}\left( x_{1}\right) \partial _{x_{1}}w_{n}\right)
,\upsilon \right) \bar{K}\left( \left( P_{\varepsilon }P_{+}w_{n},\upsilon
\right) \right) - 
\]%
\[
\bar{K}\left( \left( P_{\varepsilon }P_{+}a^{2}\left( x_{1}\right) \partial
_{x_{1}}w_{n}\right) ,\upsilon \right) \left( P_{\varepsilon
}P_{+}w_{n},\upsilon \right) +\left( \left( P_{\varepsilon }P_{+}ib\left(
x_{1}\right) \partial _{x_{1}}w_{n}\right) ,\upsilon \right) \bar{K}\left(
\left( P_{\varepsilon }P_{+}w_{n},\upsilon \right) \right) - 
\]%
\[
\bar{K}\left( P_{\varepsilon }P_{+}\left( ib\left( x_{1}\right) \partial
_{x_{1}}w_{n}\right) ,\upsilon \right) \left( P_{\varepsilon
}P_{+}w_{n},\upsilon \right) . 
\]

\bigskip By taking the imaginary part in $\left( 3.19\right) $ we get%
\[
\partial _{t}\left\vert \left( P_{\varepsilon }P_{+}w_{n},\upsilon \right)
\right\vert ^{2}+2\func{Im}\left( \Delta P_{\varepsilon }P_{+}w_{n},\upsilon
\right) \left( \bar{K}\left( P_{\varepsilon }P_{+}w_{n},\upsilon \right)
\right) + 
\]

\[
2\func{Im}\left( AP_{\varepsilon }P_{+}w_{n},\upsilon \right) \left( \bar{K}%
\left( P_{\varepsilon }P_{+}w_{n},\upsilon \right) \right) =2\func{Im}\left(
P_{\varepsilon }P_{+}\left( Vw_{n}\right) ,\upsilon \right) \left( \bar{K}%
\left( P_{\varepsilon }P_{+}w_{n},\upsilon \right) \right) + 
\]

\[
2\func{Im}\left( \left( P_{\varepsilon }P_{+}hw_{n}\right) ,\upsilon \right)
\left( \bar{K}\left( P_{\varepsilon }P_{+}w_{n},\upsilon \right) \right) +2%
\func{Im}\left( \left( P_{\varepsilon }P_{+}w_{n},\upsilon \right) \right) 
\bar{K}\left( P_{\varepsilon }P_{+}\left( \tilde{F}_{n}\right) ,\upsilon
\right) + 
\]

\begin{equation}
2\func{Im}\left( \left( P_{\varepsilon }P_{+}a^{2}\left( x_{1}\right)
\partial _{x_{1}}w_{n}\right) ,\upsilon \right) \left( \bar{K}\left(
P_{\varepsilon }P_{+}w_{n},\upsilon \right) \right) +  \tag{3.20}
\end{equation}%
\[
2\func{Im}\left( \left( P_{\varepsilon }P_{+}ib\left( x_{1}\right) \partial
_{x_{1}}w_{n}\right) ,\upsilon \right) \left( \bar{K}\left( P_{\varepsilon
}P_{+}w_{n},\upsilon \right) \right) . 
\]

\bigskip Since for all $n\in \mathbb{Z}^{+}$, $w_{n}\left( .\right) \in X$, $%
\tilde{F}_{n}\left( .,t\right) \in X$ and for a.e. $t\in \left[ 0,1\right] $
by integrating both sides of $\left( 3.20\right) $ on $R^{n}$ we get 
\[
\dint\limits_{R^{n}}\func{Im}\left( \Delta P_{\varepsilon
}P_{+}w_{n},\upsilon \right) \left( \bar{K}\left( P_{\varepsilon
}P_{+}w_{n},\upsilon \right) \right) dx=0. 
\]

It is clear to see that%
\[
\left( u,\upsilon \right) _{X}=\left( u\left( .\right) ,\upsilon \left(
.\right) _{H}\right) _{L^{2}\left( R^{n}\right) }\text{, for }u\text{, }%
\upsilon \in X. 
\]

\bigskip Then applying the Cauchy-Schwartz and Holder inequalites for a.e. $%
t\in \left[ 0,1\right] $ we obtain 
\begin{equation}
\dint\limits_{R^{n}}\func{Im}\left( P_{\varepsilon }P_{+}Vw_{n},\upsilon
\right) \left( \bar{K}\left( P_{\varepsilon }P_{+}w_{n},\upsilon \right)
\right) dx\leq C\left\Vert V\right\Vert _{B}\left\Vert w_{n}\right\Vert
_{X}^{2}\left\Vert \upsilon \right\Vert _{X}^{2},  \tag{3.21}
\end{equation}%
\begin{equation}
\dint\limits_{R^{n}}\func{Im}\left( P_{\varepsilon }P_{+}hw_{n},\upsilon
\right) \bar{K}\left( P_{\varepsilon }P_{+}w_{n},\upsilon \right) dx\leq
C\left\Vert h\right\Vert _{L^{\infty }}\left\Vert w_{n}\right\Vert
_{X}^{2}\left\Vert \upsilon \right\Vert _{X}^{2},  \tag{3.22}
\end{equation}%
\begin{equation}
\dint\limits_{R^{n}}\func{Im}\left( P_{\varepsilon }P_{+}\tilde{F}%
_{n}w_{n},\upsilon \right) \left( \bar{K}\left( P_{\varepsilon
}P_{+}w_{n},\upsilon \right) \right) dx\leq C\left\Vert \tilde{F}%
_{n}\right\Vert _{X}\left\Vert w_{n}\right\Vert _{X}^{2}\left\Vert \upsilon
\right\Vert _{X}^{2}.  \tag{3.23}
\end{equation}

Moreover, again applying the Cauchy-Schwartz and Holder inequalites due to
symmetricity of the operator $A,$ for a.e. $t\in \left[ 0,1\right] $ we get

\begin{equation}
\dint\limits_{R^{n}}\func{Im}\left( AP_{\varepsilon }P_{+}w_{n},\upsilon
\right) \left( \bar{K}\left( P_{\varepsilon }P_{+}w_{n},\upsilon \right)
\right) dx\leq C\left\Vert Aw_{n}\right\Vert _{X}\left\Vert w_{n}\right\Vert
_{X}\left\Vert \upsilon \right\Vert _{X}^{2}  \tag{3.24}
\end{equation}%
where, the constant $C$ in $\left( 3.21\right) -\left( 3.24\right) $ is
intependent of $\upsilon \in C_{0}^{\infty }\left( R^{n};H\right) ,$ $%
\varepsilon \in \left( 0,\left. 1\right] \right. $ and $n\in \mathbb{Z}^{+}.$

Since $C_{0}^{\infty }\left( R^{n};H\right) $ is dense in $X,$ from $\left(
3.24\right) $-$\left( 3.27\right) $ in view of operator theory in Hilbert
spaces, we obtain the following%
\[
\left\vert \dint\limits_{R^{n}}\func{Im}\left( P_{\varepsilon }P_{+}Vw_{n},%
\bar{K}P_{\varepsilon }P_{+}w_{n}\right) dx\right\vert \leq C\left\Vert
V\right\Vert _{B}\left\Vert w_{n}\right\Vert _{X}^{2}, 
\]%
\[
\left\vert \dint\limits_{R^{n}}\func{Im}\left( P_{\varepsilon }P_{+}hw_{n},%
\bar{K}P_{\varepsilon }P_{+}w_{n}\right) dx\right\vert \leq C\left\Vert
h\right\Vert _{L^{\infty }}\left\Vert w_{n}\right\Vert _{X}^{2}\leq C\frac{1%
}{n}\left\Vert w_{n}\right\Vert _{X}^{2}, 
\]%
\begin{equation}
\left\vert \dint\limits_{R^{n}}\func{Im}\left( P_{\varepsilon }P_{+}\tilde{F}%
_{n}w_{n},\bar{K}P_{\varepsilon }P_{+}w_{n}\right) dx\right\vert \leq
C\left\Vert \tilde{F}_{n}\right\Vert _{X}\left\Vert w_{n}\right\Vert
_{X}^{2},  \tag{3.25}
\end{equation}

\[
\left\vert \dint\limits_{R^{n}}\func{Im}\left( AP_{\varepsilon }P_{+}w_{n},%
\bar{K}P_{\varepsilon }P_{+}w_{n}\right) dx\right\vert \leq C\left\Vert
Aw_{n}\right\Vert _{X}\left\Vert w_{n}\right\Vert _{X}^{2} 
\]

For bounding the last two terms in $\left( 3.20\right) $ we will use the
abstract version of Calder\'{o}n's first commutator estimates $\left[ 4%
\right] .$ Really by Cauchy-Schvartz inequality and in view of Theorem A$%
_{1} $we get 
\begin{equation}
\left\Vert \left( \left[ P_{\pm };a\right] \partial _{x_{1}}f,\upsilon
\right) \right\Vert _{X}\leq C\left\Vert \partial _{x_{1}}a\right\Vert
_{L^{\infty }}\left\Vert f\right\Vert _{X}\left\Vert \upsilon \right\Vert
_{X},\text{ }  \tag{3.26}
\end{equation}%
\begin{equation}
\left\Vert \partial _{x_{1}}\left( \left[ P_{\pm };a\right] f,\upsilon
\right) \right\Vert _{X}\leq C\left\Vert \partial _{x_{1}}a\right\Vert
_{L^{\infty }}\left\Vert f\right\Vert _{X}\left\Vert \upsilon \right\Vert
_{X},\text{ }  \tag{3.27}
\end{equation}

\bigskip Also, from the calculus of pseudodifferential operators with
operator coefficients (see e.g. $\left[ 5\right] $ ) and the inequality $%
(3.15)$, we have

\begin{equation}
\left\Vert \left( \left[ P_{\varepsilon };a\right] \partial
_{x_{1}}f,\upsilon \right) \right\Vert _{X}\leq \frac{C}{n}\left\Vert
f\right\Vert _{X}\left\Vert \upsilon \right\Vert _{X},\text{ }  \tag{3.28}
\end{equation}%
\begin{equation}
\left\Vert \partial _{x_{1}}\left( \left[ P_{\varepsilon };a\right]
f,\upsilon \right) \right\Vert _{X}\leq \frac{C}{n}\left\Vert f\right\Vert
_{X}\left\Vert \upsilon \right\Vert _{X},\text{ }  \tag{3.29}
\end{equation}%
where the constant $C$ in $\left( 3.26\right) -(3.29)$ is independent of $%
\varepsilon \in \lbrack 0,1]$ and $n.$

We remark that estimates $(3.26)-(3.29)$ also hold with $b\left(
x_{1}\right) $ replacing $a(x_{1})$. Since $C_{0}^{\infty }\left(
R^{n};H\right) $ is dense in $X,$ from $\left( 3.26\right) $-$\left(
3.29\right) $ in view of operator theory in Hilbert spaces, we obtain

\[
\left\Vert \left[ P_{\pm };a\right] \partial _{x_{1}}f\right\Vert _{X}\leq
C\left\Vert \partial _{x_{1}}a\right\Vert _{L^{\infty }}\left\Vert
f\right\Vert _{X},\text{ }\left\Vert \partial _{x_{1}}\left[ P_{\pm };a%
\right] f\right\Vert _{X}\leq C\left\Vert \partial _{x_{1}}a\right\Vert
_{L^{\infty }}\left\Vert f\right\Vert _{X}, 
\]%
\begin{equation}
\left\Vert \left( \left[ P_{\varepsilon };a\right] \partial
_{x_{1}}f,\upsilon \right) \right\Vert _{X}\leq \frac{C}{n}\left\Vert
f\right\Vert _{X},\text{ }\left\Vert \partial _{x_{1}}\left( \left[
P_{\varepsilon };a\right] f,\upsilon \right) \right\Vert _{X}\leq \frac{C}{n}%
\left\Vert f\right\Vert _{X}\text{ ,}  \tag{3.30}
\end{equation}%
and the same estimates $\left( 3.30\right) $ with $b(x_{1}$) replacing $%
a(x_{1})$.

By reasoning as in $\left[ \text{4, Lemma 2.1}\right] $ (claim 1 and 2 )
from $\left( 3.30\right) $\ we obtain 
\begin{equation}
\left\vert \func{Im}\left( \left( P_{\varepsilon }P_{+}a^{2}\left(
x_{1}\right) \partial _{x_{1}}w_{n},\bar{K}P_{\varepsilon }P_{+}w_{n}\right)
\right) \right\vert \leq O\left( n^{-1}\left\Vert w_{n}\right\Vert
_{X}\right) ,  \tag{3.31}
\end{equation}

\[
\left\vert \func{Im}\left( \left( P_{\varepsilon }P_{+}b\left( x_{1}\right)
\partial _{x_{1}}w_{n},\bar{K}P_{\varepsilon }P_{+}w_{n}\right) \right)
\right\vert \leq O\left( n^{-1}\left\Vert w_{n}\right\Vert _{X}\right) . 
\]

Now, the estimates $\left( 3.25\right) $ and $\left( 3.31\right) $ implay
the assertion.

\textbf{4.} \textbf{Proof of Theorem 1.}

\textbf{\ }We will apply Lemma 3.1 to a solution of the equation $(3.2)$.
Since $0<\alpha <\beta =$ $\beta (k)$ for $k>k_{0}$ it follows that $\alpha
\leq \sigma \left( t\right) \leq \beta $ for any $t\in \lbrack 0,1]$.
Therefore if $y=\sqrt{\alpha \beta }x\sigma \left( t\right) $, then from $%
\left( 3.8\right) $ we get

\begin{equation}
\sqrt{\alpha \beta ^{-1}}\left\vert x\right\vert \leq \left\vert
y\right\vert \sqrt{\alpha ^{-1}\beta }\left\vert x\right\vert =\left(
ka_{0}^{-1}\right) ^{\frac{1}{2p}}\left\vert x\right\vert  \tag{4.1}
\end{equation}

\bigskip Thus, 
\begin{equation}
\left\Vert \alpha \beta \sigma ^{2}\left( t\right) V\left( \sqrt{\alpha
\beta }x\sigma \left( t\right) ,\beta t\sigma \left( t\right) \right)
\right\Vert _{L\left( H\right) }\leq \alpha ^{-1}\beta \left\Vert
V\right\Vert _{B}=\left( ka_{0}^{-1}\right) ^{\frac{1}{p}}\left\Vert
V\right\Vert _{B}  \tag{4.2}
\end{equation}%
and so, 
\begin{equation}
\left\Vert \tilde{V}\left( .,t\right) \right\Vert _{L^{\infty }\left(
R^{n};H\right) }\leq \left( ka_{0}^{-1}\right) ^{\frac{1}{p}}\left\Vert
V\left( .,t\right) \right\Vert _{L^{\infty }\left( R^{n};H\right) }. 
\tag{4.3}
\end{equation}

Also, for $s=\beta t\sigma \left( t\right) $ it is clear that 
\begin{equation}
\frac{ds}{dt}=\alpha \beta \sigma ^{2}\left( t\right) \text{, }dt=\left(
\alpha \beta \right) ^{-1}\sigma ^{-2}\left( t\right) ds.  \tag{4.4}
\end{equation}

Therefore, 
\[
\dint\limits_{0}^{1}\left\Vert \tilde{V}\left( .,t\right) \right\Vert
_{L^{\infty }\left( R^{n};H\right) }dt=\dint\limits_{0}^{1}\left\Vert
V\left( .,s\right) \right\Vert _{L^{\infty }\left( R^{n};H\right) }ds, 
\]%
and from $\left( 4.1\right) $ we get 
\begin{equation}
\dint\limits_{0}^{1}\left\Vert \tilde{V}\left( .,t\right) \right\Vert
_{L^{\infty }\left( \left\vert x\right\vert >R;H\right)
}dt=\dint\limits_{0}^{1}\left\Vert V\left( .,s\right) \right\Vert
_{L^{\infty }\left( \left\vert y\right\vert >\varkappa ;H\right) }ds, 
\tag{4.5}
\end{equation}%
where 
\[
\varkappa =\left( a_{0}k^{-1}\right) ^{\frac{1}{2p}}R. 
\]

So, if 
\[
\dint\limits_{0}^{1}\left\Vert V\left( .,s\right) \right\Vert _{L^{\infty
}\left( \left\vert y\right\vert >\varkappa ;H\right) }ds<\varepsilon _{0} 
\]%
then, 
\[
\dint\limits_{0}^{1}\left\Vert \tilde{V}\left( .,t\right) \right\Vert
_{L^{\infty }\left( \left\vert y\right\vert >R;H\right) }ds<\varepsilon _{0},%
\text{ for }R=\varkappa \left( ka_{0}^{-1}\right) ^{\frac{1}{2p}} 
\]%
and we can applay Lemma 3.1 to the equation $\left( 3.2\right) $ with%
\[
\tilde{V}\mathbb{=}\tilde{V}_{\chi \left( \left\vert x\right\vert >R\right)
}\left( x,t\right) \text{, }\tilde{F}\mathbb{=}\tilde{V}_{\chi \left(
\left\vert x\right\vert <R\right) }\left( x,t\right) \tilde{u}\left(
x,t\right) 
\]%
to get the following estimate 
\[
\sup\limits_{t\in \left[ 0,1\right] }\left\Vert e^{\nu }\tilde{u}\left(
.,t\right) \right\Vert _{X}\leq M_{0}\left( \left\Vert e^{\nu }\tilde{u}%
\left( .,0\right) \right\Vert _{X}+\left\Vert e^{\nu }\tilde{u}\left(
.,1\right) \right\Vert _{X}\right) + 
\]%
\begin{equation}
M_{0}e^{M}e^{\nu _{0}}\left\Vert \tilde{V}\right\Vert _{B}\left\Vert u\left(
.,0\right) \right\Vert _{X},  \tag{4.6}
\end{equation}%
where $M$ a positive constant defined in Remark 2.1 and 
\[
B=L^{\infty }\left( R^{n}\times \left[ 0,1\right] ;L\left( H\right) \right) 
\text{, }\nu =\left( 2p\right) ^{\frac{1}{p}}\gamma ^{\frac{1}{p}}\lambda .%
\frac{x}{2},\text{ }\nu _{0}=\left\vert \lambda \right\vert \left( 2p\right)
^{\frac{1}{p}}\gamma ^{\frac{1}{p}}\frac{R}{2}. 
\]

From $\left( 4.6\right) $ we have 
\[
\sup\limits_{t\in \left[ 0,1\right] }\dint\limits_{R^{n}}\left\Vert e^{\nu }%
\tilde{u}\left( .,t\right) \right\Vert _{H}^{2}dx\leq
M_{0}\dint\limits_{R^{n}}e^{\nu }\left( \left\Vert \tilde{u}\left(
.,0\right) \right\Vert _{H}^{2}+\left\Vert \tilde{u}\left( .,1\right)
\right\Vert _{H}^{2}\right) dx+ 
\]%
\[
M_{0}e^{M}e^{\left\vert \lambda \right\vert \left( 2p\right) ^{\frac{1}{p}%
}\gamma ^{\frac{1}{p}}R}\left\Vert \tilde{V}\right\Vert _{B}\left\Vert
u\left( .,0\right) \right\Vert _{X}^{2}, 
\]%
and multiply the above inequality by $e^{\left\vert \lambda \right\vert
/q}\left\vert \lambda \right\vert ^{n\left( q-2\right) /2}$, integrate in $%
\lambda $ and in $x$, use Fubini theorem and the following formula 
\begin{equation}
e^{\gamma \left\vert x\right\vert ^{p}/p}\thickapprox
\dint\limits_{R^{n}}e^{\gamma ^{\frac{1}{p}}\lambda .x-\left\vert \lambda
\right\vert ^{q/q}}\left\vert \lambda \right\vert ^{n\left( q-2\right)
/2}d\lambda ,  \tag{4.7}
\end{equation}%
proven in $\left[ \text{7, Appendx}\right] $ to obtain 
\[
\dint\limits_{\left\vert x\right\vert >1}e^{2\gamma \left\vert x\right\vert
^{p}}\left\Vert \tilde{u}\left( .,t\right) \right\Vert _{H}^{2}dx\leq
M_{0}\dint\limits_{R^{n}}e^{2\gamma \left\vert x\right\vert ^{p}}\left(
\left\Vert \tilde{u}\left( .,0\right) \right\Vert _{H}^{2}+\left\Vert \tilde{%
u}\left( .,1\right) \right\Vert _{H}^{2}\right) dx+ 
\]%
\begin{equation}
M_{0}e^{2M}e^{2\gamma R^{p}}R^{C_{p}}\left\Vert \tilde{V}\right\Vert
_{B}\left\Vert u\left( .,0\right) \right\Vert _{X}^{2}.  \tag{4.8}
\end{equation}

Hence, the esimates $\left( 3.6\right) $, $\left( 3.8\right) $, $\left(
3.9\right) ,$ $\left( 4.3\right) $ and $\left( 4.8\right) $ imply 
\[
\sup\limits_{t\in \left[ 0,1\right] }\left\Vert e^{\gamma \left\vert
x\right\vert ^{p}}\tilde{u}\left( .,t\right) \right\Vert _{X}\leq
M_{0}\left( \left\Vert e^{\gamma \left\vert x\right\vert ^{p}}\tilde{u}%
\left( .,0\right) \right\Vert _{X}+\left\Vert e^{\gamma \left\vert
x\right\vert ^{p}}\tilde{u}\left( .,1\right) \right\Vert _{X}\right) + 
\]%
\begin{equation}
M_{0}e^{M}e^{\gamma }\left\Vert u\left( .,0\right) \right\Vert
_{X}+M_{0}\left( ka_{0}^{-1}\right) ^{C_{p}}e^{M}e^{\varkappa ^{p}\gamma ^{%
\frac{1}{2}}\left( ka_{0}^{-1}\right) }\left\Vert u\left( .,0\right)
\right\Vert _{X}\left\Vert V\right\Vert _{B}\leq  \tag{4.9}
\end{equation}%
\[
M_{0}\left( A_{0}+A_{k}\right) +M_{0}e^{M}\left\Vert u\left( .,0\right)
\right\Vert _{X}\left( e^{\gamma }+\left( ka_{0}^{-1}\right)
^{C_{p}}\left\Vert V\right\Vert _{B}\right) e^{k\varkappa ^{p}}\leq 
\]%
\[
M_{0}A_{k}=M_{0}e^{a_{1}k^{1/2-p}}\text{ for }k>k_{0}\left( M_{0}\right) 
\text{ sufficiently large.} 
\]

Next, we shall obtain bounds for the $\nabla \tilde{u}$. Let $\tilde{\gamma}=%
\frac{\gamma }{2}$ and $\varphi $ be a strictly convex complex valued
function on compact sets of $R^{n}$, radial such that (see $[7]$)%
\[
D^{2}\varphi \geq p(p-1)|x|^{(p-2)},\text{ for }|x|\geq 1, 
\]%
\[
\varphi \geq 0,\text{ }\left\Vert \partial ^{\alpha }\varphi \right\Vert
_{L^{\infty }}\leq C\text{, }2\leq |\alpha |\leq 4,\text{ }\left\Vert
\partial ^{\alpha }\varphi \right\Vert _{L^{\infty }\left( \left\vert
x\right\vert <2\right) }\leq C\text{ for }|\alpha |\leq 4, 
\]%
\[
\varphi (x)=|x|^{p}+O(|x|)\text{, for }|x|>1. 
\]

Let us consider the equation 
\begin{equation}
\partial _{t}\upsilon =i\left( \Delta \upsilon +A\upsilon +F\left(
x,t\right) \right) ,\text{ }x\in R^{n},\text{ }t\in \left[ 0,1\right] , 
\tag{4.10}
\end{equation}%
where $F\left( x,t\right) =\tilde{V}\upsilon ,$ $A$ is a symmetric operator
in $H$ and $\tilde{V}$ is a operator in $H$ defined by $\left( 3.3\right) .$

Let 
\[
f(x,t)=e^{\tilde{\gamma}\varphi }\upsilon (x,t),\text{ }Q\left( t\right)
=\left( f(x,t),f(x,t)\right) _{H}, 
\]%
where $\upsilon $ is a solution of $\left( 4.10\right) $. Then, by reasoning
as in $\left[ \text{17, Lemma 3.3}\right] $ we have 
\begin{equation}
\partial _{t}f=Sf+Kf+i\left[ A+e^{\tilde{\gamma}\varphi }F\right] \text{, }%
\left( x,t\right) \in R^{n}\times \left[ 0,1\right] ,  \tag{4.11}
\end{equation}%
where $S$, $K$ are symmetric and skew-symmetric operator, respectively given
by%
\begin{equation}
S=-i\tilde{\gamma}\left( 2\nabla \varphi .\nabla +\Delta \varphi \right) 
\text{, }K=i\left( \Delta +A+\tilde{\gamma}^{2}\left\vert \nabla \varphi
\right\vert ^{2}\right) .  \tag{4.12}
\end{equation}%
Let 
\[
\left[ S,K\right] =SK-KS\text{.} 
\]%
A calculation shows that, 
\[
SK=\tilde{\gamma}\left( 2\nabla \varphi .\nabla +\Delta \varphi \right)
\left( \Delta +A+\tilde{\gamma}^{2}\left\vert \nabla \varphi \right\vert
^{2}\right) =\tilde{\gamma}\left( 2\nabla \varphi .\nabla +\Delta \varphi
\right) \Delta + 
\]%
\[
\tilde{\gamma}\left( 2\nabla \varphi .\nabla +\Delta \varphi \right) A+%
\tilde{\gamma}^{3}\left\vert \nabla \varphi \right\vert ^{2}\left( 2\nabla
\varphi .\nabla +\Delta \varphi \right) , 
\]%
\[
KS=\tilde{\gamma}\left[ \Delta \left( 2\nabla \varphi .\nabla +\Delta
\varphi \right) +A\left( 2\nabla \varphi .\nabla +\Delta \varphi \right) %
\right] +\tilde{\gamma}^{3}\left\vert \nabla \varphi \right\vert ^{2}\left(
2\nabla \varphi .\nabla +\Delta \varphi \right) , 
\]%
\[
\left[ S,K\right] =\tilde{\gamma}\left[ \left( 2\nabla \varphi .\nabla
+\Delta \varphi \right) \Delta -\Delta \left( 2\nabla \varphi .\nabla
+\Delta \varphi \right) \right] +2\tilde{\gamma}\left( \nabla \varphi
.\nabla A-A\nabla \varphi .\nabla \right) , 
\]%
\begin{equation}
\text{ }S_{t}+\left[ S,K\right] =-2i\tilde{\gamma}\left( \nabla \varphi
.\partial _{t}\nabla +\Delta \varphi \partial _{t}\right) +  \tag{4.13}
\end{equation}%
\[
\tilde{\gamma}\left[ \left( 2\nabla \varphi .\nabla +\Delta \varphi \right)
\Delta -\Delta \left( 2\nabla \varphi .\nabla +\Delta \varphi \right) \right]
+2\tilde{\gamma}\left( \nabla \varphi .\nabla A-A\nabla \varphi .\nabla
\right) . 
\]

By $\left[ \text{17, Lemma 3.2}\right] $ 
\[
Q^{^{\prime \prime }}\left( t\right) =2\partial _{t}\func{Re}\left( \partial
_{t}f-Sf-Kf,f\right) _{X}+2\left( S_{t}f+\left[ S,K\right] f,f\right) _{X}+ 
\]

\begin{equation}
\left\Vert \partial _{t}f-Sf+Kf\right\Vert _{X}^{2}-\left\Vert \partial
_{t}f-Sf-Kf\right\Vert _{X}^{2},  \tag{4.14}
\end{equation}%
so, 
\begin{equation}
Q^{^{\prime \prime }}\left( t\right) \geq 2\partial _{t}\func{Re}\left(
\partial _{t}f-Sf-Kf,f\right) _{X}+2\left( S_{t}f+\left[ S,K\right]
f,f\right) _{X}.  \tag{4.15}
\end{equation}

Multiplying $(4.15)$ by $t(1-t)$ and integrating in $t$ we obtain 
\[
\dint\limits_{0}^{1}t(1-t)\left( S_{t}f+\left[ S,K\right] f,f\right)
_{X}dt\leq M_{0}\left[ \sup\limits_{t\in \left[ 0,1\right] }\left\Vert e^{%
\tilde{\gamma}\varphi }\upsilon \left( .,t\right) \right\Vert
_{X}+\sup\limits_{t\in \left[ 0,1\right] }\left\Vert e^{\tilde{\gamma}%
\varphi }F\left( .,t\right) \right\Vert _{X}\right] . 
\]

This computation can be justified by parabolic regularization using the fact
that we already know the decay estimate for the solution of $\left(
4.10\right) $. Hence, combining $(3.8)$, $(4.3)$ and $(4.9$) it follows that 
\[
\tilde{\gamma}\dint\limits_{0}^{1}\dint\limits_{R^{n}}t(1-t)D^{2}\varphi
\left( x,t\right) \left( \nabla f,\nabla f\right) _{H}dxdt+\tilde{\gamma}%
^{3}\dint\limits_{0}^{1}\dint\limits_{R^{n}}t(1-t)D^{2}\varphi \left(
x,t\right) \left( \nabla f,\nabla f\right) _{H}dxdt\leq 
\]%
\begin{equation}
M_{0}\left[ \sup\limits_{t\in \left[ 0,1\right] }\left\Vert e^{\tilde{\gamma}%
\varphi }\upsilon \left( .,t\right) \right\Vert _{X}\left( 1+\left\Vert 
\tilde{V}\left( .,t\right) \right\Vert _{L^{\infty }\left( R^{n};L\left(
H\right) \right) }\right) +\tilde{\gamma}\sup\limits_{t\in \left[ 0,1\right]
}\left\Vert e^{\tilde{\gamma}\varphi }\upsilon \left( .,t\right) \right\Vert
_{X}\right] \leq  \tag{4.16}
\end{equation}%
\[
M_{0}k^{C_{p}}A_{k}. 
\]

It is clear to see that%
\[
\nabla f=\tilde{\gamma}\upsilon e^{\tilde{\gamma}\varphi }\nabla \varphi +e^{%
\tilde{\gamma}\varphi }\nabla \upsilon \text{.} 
\]

So, by using the properties of $\varphi $ we get 
\[
\text{ }\left\vert e^{2\tilde{\gamma}\varphi }D^{2}\varphi \left\vert \nabla
\varphi \right\vert ^{2}\right\vert \leq C_{p}e^{3\gamma \varphi \mid 2}. 
\]

From here, we can conclude that 
\[
\gamma \dint\limits_{0}^{1}\dint\limits_{R^{n}}t(1-t)\left( 1+\left\vert
x\right\vert \right) ^{p-2}\left\Vert \nabla \upsilon \left( x,t\right)
\right\Vert _{H}e^{\gamma \left\vert x\right\vert
^{p}}dxdt+\sup\limits_{t\in \left[ 0,1\right] }\left\Vert e^{\frac{\gamma
\left\vert x\right\vert ^{p}}{2}}\upsilon \left( .,t\right) \right\Vert
_{X}\leq 
\]%
\begin{equation}
C_{0}k^{C_{p}}A_{k}^{2}=C_{0}k^{C_{p}}e^{a\left( k,p\right) }  \tag{4.17}
\end{equation}%
for $k\geq k_{0}(M_{0})$ sufficiently large, where 
\[
a\left( k,p\right) =C_{\mu }M_{0}k^{C_{p}}e^{2a_{1}k^{1\mid \left(
2-p\right) }}. 
\]%
For proving Theorem 1 first, we deduce the following estimate%
\begin{equation}
\dint\limits_{\left\vert x\right\vert <\frac{R}{2}}\dint\limits_{\nu
_{1}}^{\nu _{2}}\left\Vert \tilde{u}\left( x,t\right) \right\Vert
_{H}dtdx\geq C_{\nu }e^{-M}\left\Vert u\left( .,0\right) \right\Vert _{X}, 
\tag{4.18}
\end{equation}%
for $R$ sufficiently large, $\nu _{1}$, $\nu _{2}\in \left( 0,1\right) $, $%
\nu _{1}<\nu _{2}$ and $\nu =\left( \nu _{1},\nu _{2}\right) $. From $\left(
3.1\right) $ by using the change of variables $s=\beta t\sigma \left(
t\right) $ and $y=\sqrt{\alpha \beta }x\sigma \left( t\right) $ we get 
\[
\dint\limits_{\left\vert x\right\vert <\frac{R}{2}}\dint\limits_{\nu
_{1}}^{\nu _{2}}\left\Vert \tilde{u}\left( x,t\right) \right\Vert
_{H}^{2}dtdx=\left( \alpha \beta \right) ^{\frac{n}{2}}\dint\limits_{\left%
\vert x\right\vert <\frac{R}{2}}\dint\limits_{\nu _{1}}^{\nu _{2}}\left\vert
\sigma \left( t\right) \right\vert ^{n}\left\Vert u\left( \sqrt{\alpha \beta 
}x\sigma \left( t\right) ,\beta t\sigma \left( t\right) \right) \right\Vert
_{H}^{2}dtdx\geq 
\]%
\begin{equation}
M_{0}\frac{\beta }{\alpha }\dint\limits_{\left\vert y\right\vert
<R_{0}}\dint\limits_{s\nu _{1}}^{s\nu _{2}}\left\Vert u\left( y,s\right)
\right\Vert _{H}^{2}\frac{dsdy}{s^{2}}\geq M_{0}\frac{\beta }{\alpha }%
\dint\limits_{\left\vert y\right\vert <R_{0}}\dint\limits_{s\nu _{1}}^{s\nu
_{2}}\left\Vert u\left( y,s\right) \right\Vert _{H}^{2}dsdy  \tag{4.19}
\end{equation}%
for $k>M_{0},$ $s\nu _{1}>\frac{1}{2}$ and $R_{0}=R\left( ka_{0}^{-1}\right)
^{\frac{1}{2p}}.$ Thus, taking 
\begin{equation}
R>\omega \left( ka_{0}^{-1}\right) ^{\frac{1}{2p}}  \tag{4.20}
\end{equation}%
with $\omega =\omega (u)$ a constant to be determined, it follows that%
\[
\Phi \geq M_{0}\frac{\beta }{\alpha }\dint\limits_{\left\vert y\right\vert
<\omega }\dint\limits_{s\nu _{1}}^{s\nu _{2}}\left\Vert u\left( y,s\right)
\right\Vert _{H}^{2}dsdy, 
\]%
where the interval $I=I_{k}=[s\nu _{1},s\nu _{2}]$ satisfies $I\subset
\lbrack 1/2,1]$ for $k$ sufficiently large. Moreover, given $\varepsilon >0$
there exists $k_{0}(\varepsilon )>0$ such that for any $k\geq k_{0}$ one has
that $I_{k}\subset \lbrack 1-\varepsilon ,1]$. By hypothesis on $u(x,t)$,
i.e. the continuity of $\left\Vert u(\text{\textperiodcentered }%
,s)\right\Vert _{X}$ at $s=1$, it follows that there exists $\omega >1$ and $%
K_{0}=K_{0}(u)$ such that for any $k\geq K_{0}$ and for any $s\in I_{k}$ 
\[
\dint\limits_{\left\vert y\right\vert <\omega }\left\Vert u\left( y,s\right)
\right\Vert _{H}^{2}dy\geq C_{\nu }e^{-M}\left\Vert u\left( .,0\right)
\right\Vert _{X}, 
\]%
which yields the desired result. Next, we deduce the following estimate%
\begin{equation}
\dint\limits_{\left\vert x\right\vert <R}\dint\limits_{\mu _{1}}^{\mu
_{2}}\left( \left\Vert \tilde{u}\left( x,t\right) \right\Vert
_{H}^{2}+\left\Vert \nabla \tilde{u}\left( x,t\right) \right\Vert
_{H}^{2}+\left\Vert A\tilde{u}\left( x,t\right) \right\Vert _{H}^{2}\right)
dtdx\leq C_{n\mu }C_{0}k^{C_{p}}e^{a\left( k,p\right) },  \tag{4.21}
\end{equation}%
for $R$ sufficiently large, $\nu _{1},$ $\nu _{2}\in \left( 0,1\right) $, $%
\mu _{1}=\frac{\left( \nu _{2}-\nu _{1}\right) }{8}$, $\mu _{2}=1-\mu _{1},$ 
$\mu _{1}<\mu _{2}$ and $\mu =\left( \mu _{1},\mu _{2}\right) $.

Indeed, from $\left( 3.9\right) $ and $\left( 4.17\right) $\ we obtain 
\[
\dint\limits_{\left\vert x\right\vert <R}\dint\limits_{\mu _{1}}^{\mu
_{2}}\left\Vert \tilde{u}\left( x,t\right) \right\Vert _{H}dtdx\leq C_{\mu
}e^{2M}\left\Vert u\left( .,0\right) \right\Vert _{X}, 
\]%
\begin{equation}
\dint\limits_{\mu _{1}}^{\mu _{2}}\dint\limits_{\left\vert x\right\vert
<R}\left\Vert \nabla \tilde{u}\left( x,t\right) \right\Vert _{H}dtdx\leq
C_{n\mu }\dint\limits_{\mu _{1}}^{\mu _{2}}\dint\limits_{\left\vert
x\right\vert <R}t(1-t)\left\Vert \nabla \upsilon \left( x,t\right)
\right\Vert _{H}e^{\gamma \left\vert x\right\vert ^{p}}dtdx\leq  \tag{4.22}
\end{equation}%
\[
C_{\mu }\gamma ^{-1}R^{2-p}C_{0}k^{C_{p}}A_{k}^{2}\leq C_{\mu
}C_{0}k^{C_{p}}e^{2a_{1}k^{\left( 2-p\right) ^{-1}}}. 
\]

Hence, from $\left( 4.22\right) $ we get $\left( 4.21\right) $ for $k\geq
k_{0}(C_{0})$ sufficiently large.

Let $Y=L^{2}\left( R^{n}\times \left[ 0,1\right] ;H\right) $. By reasoning
as in $\left[ 6\text{, Lemma 3.1}\right] $ we obtain

\textbf{Lemma 4.1}. Let $A$ be a pozitive operator in the Hilbert space $H$
and $iA$ generates a semigrop $U\left( t\right) =e^{iAt}$. Assume that $R>0$
and $\varphi $ : $[0,1]\rightarrow \mathbb{R}$ is a smooth function. Then,
there exists $C=C(n,\varphi ,H,A)>0$ such that, the inequality%
\[
\frac{\varkappa ^{\frac{3}{2}}}{R^{2}}\left\Vert e^{\varkappa \left\vert
\psi \right\vert }g\right\Vert _{Y}\leq C\left\Vert e^{\varkappa \left\vert
\psi \right\vert }i\left( \partial _{t}g+\Delta g+Ag\right) \right\Vert _{Y} 
\]%
holds, for $\varkappa \geq CR^{2}$ and $g\in C_{0}^{\infty }\left(
R^{n+1};H\right) $ with support contained in the set 
\[
\left\{ x,t:\text{ }\left\vert \psi \left( x,t\right) \right\vert
=\left\vert \frac{x}{R}+\varphi \left( t\right) e_{1}\right\vert \geq
1\right\} . 
\]

\textbf{Proof. }Let $f=e^{\alpha \left\vert \psi \left( x,t\right)
\right\vert ^{2}}g.$ Then, by acts of Schredinger operator $\left( i\partial
_{t}+\Delta +A\right) $ to $f\in X$ we get 
\[
e^{\alpha \left\vert \psi \left( x,t\right) \right\vert ^{2}}\left( \left(
i\partial _{t}g+\Delta g+Ag\right) \right) =S_{\alpha }f-4\alpha A_{\alpha
}f, 
\]

where 
\[
S_{\alpha }=\left( i\partial _{t}+\Delta +A\right) +\frac{4\alpha ^{2}}{R^{2}%
}\text{ }\left\vert \psi \left( x,t\right) \right\vert ^{2},\text{ } 
\]%
\[
A_{\alpha }=\frac{1}{R}\psi \left( x,t\right) .\nabla +\frac{n}{R^{2}}+\frac{%
i\varphi ^{\prime }}{2}\left( \frac{x}{R}+\varphi \left( t\right) \right) . 
\]%
Hence, 
\[
\left( S_{\alpha }\right) ^{\ast }=S_{\alpha }\text{, }\left( A_{\alpha
}\right) ^{\ast }=A_{\alpha } 
\]%
and 
\[
\left\Vert e^{\alpha \left\vert \psi \left( x,t\right) \right\vert
^{2}}\left( i\partial _{t}g+\Delta g+Ag\right) \right\Vert _{X}^{2}=\left(
S_{\alpha }f-4\alpha A_{\alpha }f,S_{\alpha }f-4\alpha A_{\alpha }f\right)
_{X}\geq 
\]%
\[
-4\alpha \left( \left( S_{\alpha }A_{\alpha }-A_{\alpha }S_{\alpha }\right)
,f\right) _{X}=-4\alpha \left( \left[ S_{\alpha },A_{\alpha }\right]
f,f\right) _{X}. 
\]

\bigskip A calculation shows that 
\[
\left[ S_{\alpha },A_{\alpha }\right] =\frac{2}{R^{2}}\Delta -\frac{4\alpha
^{2}}{R^{4}}\left\vert \frac{x}{R}+\varphi e_{1}\right\vert ^{2}-\frac{1}{2}+%
\frac{2i\varphi ^{\prime }}{R}\partial _{x_{1}} 
\]%
and 
\begin{equation}
\left\Vert e^{\alpha \left\vert \psi \left( x,t\right) \right\vert
^{2}}\left( i\partial _{t}g+\Delta g+Ag\right) \right\Vert _{X}^{2}\geq 
\tag{4.23}
\end{equation}%
\[
\frac{16\alpha ^{3}}{R^{4}}\dint \left\vert \frac{x}{R}+\varphi
e_{1}\right\vert ^{2}\left\Vert f\left( x,t\right) \right\Vert _{H}^{2}dxdt+%
\frac{8\alpha }{R^{2}}\dint \left\Vert \nabla f\left( x,t\right) \right\Vert
_{H}^{2}dxdt+ 
\]%
\[
2\alpha \dint \left[ \left( \frac{x}{R}+\varphi \right) \varphi ^{\prime
\prime }+\left( \varphi ^{\prime }\right) ^{2}\right] \left\Vert f\left(
x,t\right) \right\Vert _{H}^{2}dxdt-\frac{8\alpha i}{R}\dint \varphi
^{\prime }\partial _{x_{1}}f\bar{f}dxdt. 
\]

Hence, using the hypothesis on the support on $g$ and the Cauchy--Schwarz
inequality, the absolute value of the last two terms in $(4.23)$ can be
bounded by a fraction of the first two terms on the right-hand side of $%
(4.23)$, when

$\alpha >CR^{2}$ for some large $C$ depending on $\left\Vert \varphi
^{\prime }\right\Vert _{\infty }+\left\Vert \varphi ^{\prime \prime
}\right\Vert _{\infty }.$ This yields the assertion.

Now, from $\left( 3.3\right) $ we have 
\[
\left\Vert \tilde{V}\left( x,t\right) \right\Vert _{H}\leq \frac{\alpha }{%
\beta }\mu _{1}^{-2}\left\Vert V\right\Vert _{B}\leq \mu _{1}^{-2}a_{0}^{%
\frac{1}{p}}k^{-\frac{1}{p}}\left\Vert V\right\Vert _{B}. 
\]

Then from $\left( 4.20\right) $ we get 
\begin{equation}
\left\Vert \tilde{V}\right\Vert _{L^{\infty }\left( R^{n}\times \left[ \mu
_{1},\mu _{2}\right] ;L\left( H\right) \right) }<R.  \tag{4.24}
\end{equation}

Define 
\begin{equation}
\delta ^{2}\left( R\right) =\dint\limits_{\mu _{1}}^{\mu
_{2}}\dint\limits_{R-1\leq \left\vert x\right\vert \leq R}\left( \left\Vert 
\tilde{u}\left( x,t\right) \right\Vert _{H}^{2}+\left\Vert \nabla \tilde{u}%
\left( x,t\right) \right\Vert _{H}^{2}+\left\Vert A\tilde{u}\left(
x,t\right) \right\Vert _{H}^{2}\right) dtdx.  \tag{4.25}
\end{equation}

Let $\nu _{1}$, $\nu _{2}\in \left( 0,1\right) $, $\nu _{1}<\nu _{2}$ and $%
\nu _{2}<2\nu _{1}$. We choose $\varphi \in C^{\infty }\left( \left[ 0,1%
\right] \right) $ and $\theta $, $\theta _{R}\in C_{0}^{\infty }\left(
R^{n}\right) $ satisfying%
\[
0\leq \varphi \left( t\right) \leq 3\text{, }\varphi \left( t\right) =3\text{
for }t\in \left[ \nu _{1},\nu _{2}\right] \text{, }\varphi \left( t\right) =0%
\text{ for } 
\]%
\[
t\in \left[ 0,\nu _{2}-\nu _{1}\right] \cup \left[ \nu _{2}+\frac{\nu
_{2}-\nu _{1}}{2},1\right] , 
\]

\[
\theta _{R}\left( x\right) =1\text{ for }\left\vert x\right\vert <R-1,\text{ 
}\theta _{R}\left( x\right) =0\text{, for }\left\vert x\right\vert >R\text{,}
\]%
and 
\[
\theta \left( x\right) =1\text{ for }\left\vert x\right\vert <1,\text{ }%
\theta \left( x\right) =0\text{, for }\left\vert x\right\vert \geq 2\text{.} 
\]

Let 
\begin{equation}
g\left( x,t\right) =\theta _{R}\left( x\right) \theta \left( \psi \left(
x,t\right) \right) \tilde{u}\left( x,t\right) ,  \tag{4.26}
\end{equation}%
where $\tilde{u}\left( x,t\right) $ is a solution of $\left( 3.2\right) $
when $\tilde{V}=\tilde{F}\equiv 0$. It is clear to see that 
\[
\left\vert \psi \left( x,t\right) \right\vert \geq \frac{5}{2}\text{ for }%
\left\vert x\right\vert <\frac{R}{2}\text{ and }t\in \left[ \nu _{1},\nu _{2}%
\right] . 
\]%
Hence, 
\[
g\left( x,t\right) =\tilde{u}\left( x,t\right) \text{ and }e^{\varkappa
\left\vert \psi \left( x,t\right) \right\vert ^{2}}\geq e^{\frac{25}{4}%
\varkappa }\text{ for }\left\vert x\right\vert <\frac{R}{2}\text{, }t\in %
\left[ \nu _{1},\nu _{2}\right] . 
\]%
Moreover, from $\left( 4.26\right) $ also we get that 
\[
g\left( x,t\right) =0\text{ for }\left\vert x\right\vert \geq R\text{ or }%
t\in \left[ 0,\nu _{2}-\nu _{1}\right] \cup \left[ \nu _{2}+\frac{\nu
_{2}-\nu _{1}}{2},1\right] , 
\]%
so 
\[
supp\text{ }g\subset \left\{ \left\vert x\right\vert \leq R\right\} \times %
\left[ \nu _{2}-\nu _{1},\nu _{2}+\frac{\nu _{2}-\nu _{1}}{2}\right] \cap
\left\{ \left\vert b\left( x,t\right) \right\vert \geq 1\right\} . 
\]

Then, for $\xi =\psi \left( x,t\right) $ we have%
\[
\left( i\partial _{t}+\Delta +A+\tilde{V}\right) g=\left[ \theta \left( \xi
\right) \left( 2\nabla \theta _{R}\left( x\right) .\tilde{u}+\tilde{u}\Delta
\theta _{R}\left( x\right) \right) +2\nabla \theta \left( \xi \right)
.\nabla \theta _{R}\tilde{u}\right] + 
\]%
\[
\theta _{R}\left( x\right) \left[ 2R^{-1}\nabla \theta \left( \xi \right)
.\nabla \tilde{u}+R^{-2}\tilde{u}\Delta \theta \left( \xi \right) +i\varphi
^{\prime }\partial x_{1}\theta \left( \xi \right) u\right] =B_{1}+B_{2}\text{%
.} 
\]

Note that, 
\[
supp\text{ }B_{1}\subset \left\{ \left( x,t\right) :\text{ }R-1\leq
\left\vert x\right\vert \leq R,\text{ }\mu _{1}\leq t\leq \mu _{2}\right\} 
\]

and 
\[
supp\text{ }B_{2}\subset \left\{ \left( x,t\right) \in R^{n}\times \left[ 0,1%
\right] \text{, }1\leq \left\vert \psi \left( x,t\right) \right\vert \leq
2\right\} . 
\]

Now applying Lemma 4.1 choosing $\varkappa =d_{n}^{2}R^{2},$ $d_{n}^{2}\geq
\left\Vert \varphi ^{\prime }\right\Vert _{\infty }+\left\Vert \varphi
^{\prime \prime }\right\Vert _{\infty }$ it follows that

\begin{equation}
R\left\Vert e^{\varkappa \left\vert \psi \right\vert ^{2}}g\right\Vert
_{Y}\leq C\left\Vert e^{\varkappa \left\vert \psi \right\vert ^{2}}i\left(
\partial _{t}g+\Delta g+Ag\right) \right\Vert _{Y}\leq  \tag{4.27}
\end{equation}

\[
C\left[ \left\Vert e^{\varkappa \left\vert \psi \right\vert ^{2}}\tilde{V}%
g\right\Vert _{Y}+\left\Vert e^{\varkappa \left\vert \psi \right\vert
^{2}}B_{1}\right\Vert _{Y}+\left\Vert e^{\varkappa \left\vert \psi
\right\vert ^{2}}B_{2}\right\Vert _{Y}\right] = 
\]%
\[
D_{1}+D_{2}+D_{3}. 
\]

Since 
\[
\left\Vert \tilde{V}\right\Vert _{L^{\infty }\left( R^{n}\times \left[ \mu
_{1},\mu _{2}\right] ;L\left( H\right) \right) }<R, 
\]%
$D_{1}$ can be absorbed in the left hand side of $\left( 4.27\right) .$
Moreover, $\left\vert \psi \left( x,t\right) \right\vert \leq 4$ on the
support of $B_{1},$ thus \ 
\[
D_{2}\leq C\delta \left( R\right) e^{16\varkappa }. 
\]

Let%
\[
R_{\mu }^{n}=\left\{ \left( x,t\right) :\left\vert x\right\vert \leq R,\text{
}\mu _{1}\leq t\leq \mu _{2}\right\} 
\]%
Then $R_{\mu }^{n}\subset $supp $B_{2}$, and $1\leq \left\vert \psi \left(
x,t\right) \right\vert \leq 2$, so%
\[
D_{3}\leq Ce^{4\varkappa }\left\Vert \tilde{u}+\nabla \tilde{u}\right\Vert
_{L^{2}\left( R^{n}\times \left[ \mu _{1},\mu _{2}\right] ;H\right) }. 
\]%
By using $\left( 4.18\right) $ and $\left( 4.22\right) $ we have 
\begin{equation}
C_{\mu }e^{-M}e^{\frac{25}{4}\varkappa }\left\Vert u\left( .,0\right)
\right\Vert _{X}\leq Re^{\frac{25}{4}\varkappa }\left[ \dint\limits_{\left%
\vert x\right\vert <\frac{R}{2}}\dint\limits_{\nu _{1}}^{\nu _{2}}\left\Vert 
\tilde{u}\left( x,t\right) \right\Vert _{H}dtdx\right] ^{\frac{1}{2}}\leq 
\tag{4.28}
\end{equation}%
\[
C_{\mu }\delta \left( R\right) e^{\frac{25}{4}\varkappa }+C_{\mu
}e^{4\varkappa }\left\Vert \tilde{u}+\nabla \tilde{u}\right\Vert
_{L^{2}\left( R_{\mu }^{n};H\right) }\leq C_{\mu }\delta \left( R\right)
e^{16\varkappa }+C_{\mu }C_{0}k^{C_{p}}e^{4\varkappa }e^{2a_{1}k^{\left(
2-p\right) ^{-1}}}. 
\]

\bigskip Puting $\varkappa =d_{n}R^{2}=2a_{1}k^{\frac{1}{2-p}}$ it follows
from $\left( 4.28\right) $ that, if $\left\Vert u\left( .,0\right)
\right\Vert _{X}\neq 0$ then%
\begin{equation}
\delta \left( R\right) \geq C_{\mu }\left\Vert u\left( .,0\right)
\right\Vert _{X}e^{-\left( M+10\varkappa \right) }=C_{\mu }\left\Vert
u\left( .,0\right) \right\Vert _{X}e^{-\left( M+20\right) a_{1}k^{\frac{1}{%
2-p}}}  \tag{4.29 }
\end{equation}%
for $k\geq k_{0}(C_{\mu })$ sufficiently large.

Now, by $\left( 4.22\right) $ we get%
\[
\delta ^{2}\left( R\right) =\dint\limits_{\mu _{1}}^{\mu
_{2}}\dint\limits_{R-1\leq \left\vert x\right\vert \leq R}\left( \left\Vert 
\tilde{u}\left( x,t\right) \right\Vert _{H}^{2}+\left\Vert \nabla \tilde{u}%
\left( x,t\right) \right\Vert _{H}^{2}+\left\Vert A\tilde{u}\left(
x,t\right) \right\Vert _{H}^{2}\right) dtdx\leq 
\]%
\[
\dint\limits_{\mu _{1}}^{\mu _{2}}\dint\limits_{R-1\leq \left\vert
x\right\vert \leq R}\left( \left\Vert \tilde{u}\left( x,t\right) \right\Vert
_{H}^{2}+\left\Vert A\tilde{u}\left( x,t\right) \right\Vert _{H}^{2}\right)
dtdx+ 
\]%
\[
C_{\mu }\dint\limits_{\mu _{1}}^{\mu _{2}}\dint\limits_{R-1\leq \left\vert
x\right\vert \leq R}t\left( 1-t\right) \left\Vert \nabla \tilde{u}\left(
x,t\right) \right\Vert _{H}^{2}dtdx\leq 
\]%
\[
C_{\mu }e^{-\gamma \left( R-1\right) ^{p}}\sup\limits_{t\in \left[ 0,1\right]
}\left\Vert e^{\gamma \left\vert x\right\vert ^{p/2}}\tilde{u}\left(
x,t\right) \right\Vert _{X}^{2}+C_{\mu }\gamma ^{-1}R^{2-p}e^{-\gamma \left(
R-1\right) ^{p}}\times 
\]%
\begin{equation}
\dint\limits_{\mu _{1}}^{\mu _{2}}\dint\limits_{R-1\leq \left\vert
x\right\vert \leq R}\left[ \frac{t\left( 1-t\right) }{\left( 1+\left\vert
x\right\vert \right) ^{2/p}}\left\Vert \nabla \tilde{u}\left( x,t\right)
\right\Vert _{H}^{2}+\left\Vert A\tilde{u}\left( x,t\right) \right\Vert
_{H}^{2}\right] dtdx\leq  \tag{4.30}
\end{equation}%
\[
C_{\mu }\gamma ^{-1}k^{C_{p}}e^{\eta \left( p\right) },\text{ }\eta \left(
p\right) =2a_{1}k^{\left( 2-p\right) ^{-1}}-\gamma \left( R-1\right) ^{p}. 
\]

The estimates $\left( 4.28\right) -\left( 4.30\right) $ imply 
\begin{equation}
C_{\mu }e^{-2M}e^{\frac{25}{4}\varkappa }\left\Vert u\left( .,0\right)
\right\Vert _{X}\leq C_{0}k^{C_{p}}e^{\omega \left( p\right) }+O\left(
k^{1/2\left( 2-p\right) }\right) \text{,}  \tag{4.31}
\end{equation}

where 
\[
\omega \left( p\right) =42a_{1}k^{1/\left( 2-p\right) }-a_{0}^{-\frac{1}{2}%
}\left( \frac{2a_{1}}{d_{n}}\right) ^{\frac{p}{2}}k^{1/\left( 2-p\right) }. 
\]

Hence, if $42a_{1}<\sqrt{a_{1}^{p}a_{0}}\left( \frac{2}{d_{n}}\right) ^{%
\frac{p}{2}}$ by letting $k$ tends to infinity it follows from $\left(
4.31\right) $ that $\left\Vert u\left( .,0\right) \right\Vert _{X}=0$, which
gives $u\left( x,t\right) \equiv 0$.

\textbf{Proof of Corollary 1. }Since 
\[
\dint\limits_{R^{n}}\left\Vert u\left( x,1\right) \right\Vert
_{H}^{2}e^{2b\left\vert x\right\vert ^{q}}dx<\infty \text{ for }b=\frac{%
\beta ^{q}}{q} 
\]

one has that%
\[
\dint\limits_{R^{n}}\left\Vert u\left( x,1\right) \right\Vert
_{H}^{2}e^{2k\left\vert x\right\vert ^{q}}dx\leq \left\Vert e^{2k\left\vert
x\right\vert ^{q}-2b\left\vert x\right\vert ^{q}}\right\Vert _{\infty
}\dint\limits_{R^{n}}\left\Vert u\left( x,1\right) \right\Vert
_{H}^{2}e^{2b\left\vert x\right\vert ^{q}}dx. 
\]

Then, by reasoning as in $\left[ \text{7, Corollary 1}\right] $ we obtain
the assertion.

\textbf{\ Proof of Theorem 2. }Indeed, just applying Corollary 1 with 
\[
u\left( x,t\right) =u_{1}\left( x,t\right) -u_{2}\left( x,t\right) 
\]%
and 
\[
V\left( x,t\right) =\frac{F\left( u_{1},\bar{u}_{1}\right) -F\left( u_{2},%
\bar{u}_{2}\right) }{u_{1}-u_{2}} 
\]%
we obtain the assertion of Theorem 2.

\begin{center}
\textbf{5. Proof of Theorem 3. }
\end{center}

First, we deduce the corresponding upper bounds. Assume 
\[
\left\Vert u\left( .,t\right) \right\Vert _{X}=a\neq 0. 
\]

Fix $\bar{t}$ near $1,$ and let 
\[
\upsilon \left( x,t\right) =u\left( x,t-1+\bar{t}\right) \text{, }t\in \left[
0,1\right] 
\]%
which satisfies the equation $\left( 2.12\right) $ with%
\begin{equation}
\left\vert \upsilon \left( x,0\right) \right\vert \leq \frac{b_{1}}{\left( 2-%
\bar{t}\right) ^{n/2}}e^{-\frac{b_{2}\left\vert x\right\vert ^{p}}{\left( 2-%
\bar{t}\right) ^{p}}},\text{ }\left\vert \upsilon \left( x,1\right)
\right\vert \leq \frac{b_{1}}{\left( 1-\bar{t}\right) ^{n/2}}e^{-\frac{%
b_{2}\left\vert x\right\vert ^{p}}{\left( 1-\bar{t}\right) ^{p}}}  \tag{5.1}
\end{equation}%
where $A$ is a linear operator$,$ $V\left( x,t\right) $ is a given potential
operator function in a Hilbert space $H.$

From $\left( 5.1\right) $ we get%
\[
\dint\limits_{R^{n}}\left\Vert \upsilon \left( x,0\right) \right\Vert
_{H}^{2}e^{A_{0}\left\vert x\right\vert ^{q}}dx=a_{0}^{2},\text{ }%
\dint\limits_{R^{n}}\left\Vert \upsilon \left( x,1\right) \right\Vert
_{H}^{2}e^{A_{1}\left\vert x\right\vert ^{q}}dx=a_{1}^{2}, 
\]%
where 
\begin{equation}
A_{0}=\frac{b_{2}}{\left( 2-\bar{t}\right) ^{p}}\text{, }A_{1}=\frac{b_{2}}{%
\left( 1-\bar{t}\right) ^{p}}.  \tag{5.2}
\end{equation}

For $V\left( x,t\right) =F\left( u,\bar{u}\right) $, by hypothesis 
\[
\left\Vert V\left( x,t\right) \right\Vert _{H}\leq C\left\Vert u\left( x,t-1+%
\bar{t}\right) \right\Vert _{H}^{\theta }\leq \frac{C}{\left( 2-t-\bar{t}%
\right) ^{\theta n/2}}e^{\frac{C\left\vert x\right\vert ^{p}}{\left( 2-t-%
\bar{t}\right) ^{p}}}. 
\]

By using Appell transformation if we suppose that $\upsilon \left(
y,s\right) $ is a solution of 
\[
\partial _{s}\upsilon =i\left( \Delta \upsilon +A\upsilon +V\left(
y,s\right) \upsilon \right) \text{, }y\in R^{n},\text{ }t\in \left[ 0,1%
\right] , 
\]%
$\alpha $ and $\beta $ are positive, then 
\begin{equation}
\tilde{u}\left( x,t\right) =\left( \sqrt{\alpha \beta }\sigma \left(
t\right) \right) ^{\frac{n}{2}}u\left( \sqrt{\alpha \beta }x\sigma \left(
t\right) ,\beta t\sigma \left( t\right) \right) e^{\eta }.  \tag{5.3}
\end{equation}%
\ verifies the equation 
\begin{equation}
\partial _{t}\tilde{u}=i\left[ \Delta \tilde{u}+A\tilde{u}+\tilde{V}\left(
x,t\right) \tilde{u}+\tilde{F}\left( x,t\right) \right] ,\text{ }x\in R^{n},%
\text{ }t\in \left[ 0,1\right]  \tag{5.4}
\end{equation}%
with $\tilde{V}\left( x,t\right) ,$ $\tilde{F}\left( x,t\right) $ defined by 
$\left( 3.3\right) $, $\left( 3.4\right) $ and 
\begin{equation}
\left\Vert e^{\gamma \left\vert x\right\vert ^{p}}\tilde{u}_{k}\left(
x,0\right) \right\Vert _{X}=\left\Vert e^{\gamma \left( \frac{\alpha }{\beta 
}\right) ^{p/2}\left\vert x\right\vert ^{p}}\upsilon \left( x,0\right)
\right\Vert _{X}=a_{0},  \tag{5.5}
\end{equation}

\[
\left\Vert e^{\gamma \left\vert x\right\vert ^{p}}\tilde{u}_{k}\left(
x,1\right) \right\Vert _{X}=\left\Vert e^{\gamma \left( \frac{\beta }{\alpha 
}\right) ^{p/2}\left\vert x\right\vert ^{p}}\upsilon \left( x,1\right)
\right\Vert _{X}=a_{1}. 
\]%
\ It follows from expressions $\left( 3.6\right) $ and $\left( 3.8\right) $\
that 
\begin{equation}
\gamma \sim \frac{1}{\left( 1-\bar{t}\right) ^{p/2}}\text{, }\beta \sim 
\frac{1}{\left( 1-\bar{t}\right) ^{p/2}}\text{, }\alpha \sim 1.  \tag{5.6}
\end{equation}%
Next, we shall estimate 
\[
\left\Vert \tilde{V}\left( x,t\right) \right\Vert _{L_{t}^{1}L_{x}^{\infty
}\left( L\left( H\right) ,R\right) },\text{ } 
\]%
where 
\[
\text{ }L_{t}^{1}L_{x}^{\infty }\left( L\left( H\right) ,R\right)
=L^{1}\left( 0,1;L^{\infty }\left( R^{n}/O_{R}\right) ;L\left( H\right)
\right) . 
\]%
Thus, 
\[
\left\Vert \tilde{V}\left( x,t\right) \right\Vert _{L\left( H\right) }\leq 
\frac{\beta }{\alpha }\left\Vert V\left( y,s\right) \right\Vert _{L\left(
H\right) }\leq \frac{\beta }{\alpha }\frac{C}{\left( 1-\bar{t}\right)
^{\theta n/2}}e^{C\left\vert y\right\vert ^{p}}, 
\]%
with 
\[
\left\vert y\right\vert =\sqrt{\alpha \beta }\left\vert x\right\vert \sigma
\left( t\right) \geq R\sqrt{\frac{\alpha }{\beta }\sim }\frac{R}{\sqrt{\beta 
}}=CR\left( 1-\bar{t}\right) ^{1/2}. 
\]%
Hence, 
\[
\left\Vert \tilde{V}\left( .,t\right) \right\Vert _{L^{\infty }\left(
R^{n};L\left( H\right) \right) }\leq \frac{\beta }{\alpha }\left\Vert
V\left( .,s\right) \right\Vert _{L^{\infty }\left( R^{n};L\left( H\right)
\right) }\leq \frac{C}{\left( 1-\bar{t}\right) ^{1+\theta n/2}} 
\]%
and 
\begin{equation}
\left\Vert \tilde{V}\left( x,t\right) \right\Vert _{L_{t}^{1}L_{x}^{\infty
}\left( L\left( H\right) ,R\right) }\leq \frac{\beta }{\alpha }\left\Vert
V\left( y,s\right) \right\Vert _{L_{t}^{1}L_{x}^{\infty }\left( R;\beta
\right) }\leq  \tag{5.7}
\end{equation}%
\[
\frac{C}{\left( 1-\bar{t}\right) ^{1+\theta n/2}}e^{-CR^{p}\left( 1-\bar{t}%
\right) ^{1/2}}, 
\]%
where 
\[
L_{t}^{1}L_{x}^{\infty }\left( R,\beta \right) =L^{1}\left( 0,1;L^{\infty
}\left( R^{n}/O_{CR/\sqrt{\beta }}\right) ;L\left( H\right) \right) . 
\]

To apply Lemma 3.1 we need 
\begin{equation}
\left\Vert \tilde{V}\left( x,t\right) \right\Vert _{L_{t}^{1}L_{x}^{\infty
}\left( L\left( H\right) ,R\right) }\leq \frac{C}{\left( 1-\bar{t}\right)
^{1+\theta n/2}}e^{-CR^{p}\left( 1-\bar{t}\right) ^{1/2}}\leq \delta _{0}. 
\tag{5.8}
\end{equation}%
for some $R,$ i.e., 
\begin{equation}
R\sim \frac{C}{\left( 1-\bar{t}\right) ^{p/2}}\delta \left( t\right) \text{, 
}  \tag{5.9}
\end{equation}%
where 
\[
\delta \left( t\right) =\log ^{\frac{1}{p}}\phi \left( t\right) \text{, }%
\phi \left( t\right) =\frac{C}{\delta _{0}\left( 1-\bar{t}\right) ^{\theta
n/2}}. 
\]%
Let 
\begin{equation}
\mathbb{V=}\tilde{V}_{\chi \left( \left\vert x>R\right\vert \right) }\left(
x,t\right) ,\text{ }\mathbb{F}=\tilde{V}_{\chi \left( \left\vert
x<R\right\vert \right) }\left( x,t\right) \tilde{u}\left( x,t\right) . 
\tag{5.10}
\end{equation}

By using $\left( 5.5\right) -\left( 5.10\right) ,$ by virtue of Lemma 3.1
and $\left( 4.7\right) $ we deduced 
\[
\sup\limits_{t\in \left[ 0,1\right] }\left\Vert e^{\gamma \left\vert
x\right\vert ^{p}}\tilde{u}\left( .,t\right) \right\Vert _{X}^{2}\leq
C\left( \left\Vert e^{\gamma \left\vert x\right\vert ^{p}}\tilde{u}\left(
.,0\right) \right\Vert _{X}^{2}+\left\Vert e^{\gamma \left\vert x\right\vert
^{p}}\tilde{u}\left( .,1\right) \right\Vert _{X}^{2}\right) + 
\]%
\[
Ca^{2}e^{C\gamma R^{p}}\left\Vert \tilde{V}\left( x,t\right) \right\Vert
\leq \frac{Ca^{2}}{\left( 1-\bar{t}\right) ^{p}}e^{\delta \left( t\right) }, 
\]

where 
\[
a=\left\Vert u_{0}\right\Vert _{X}. 
\]

Next, using the same argument given in section $4$, $(4.10)-(4.22)$, one
finds that 
\[
\gamma \dint\limits_{0}^{1}\dint\limits_{R^{n}}t(1-t)\left( 1+\left\vert
x\right\vert \right) ^{p-2}\left\Vert \nabla \tilde{u}\left( x,t\right)
\right\Vert _{H}e^{\gamma \left\vert x\right\vert ^{p}}dxdt\leq \frac{Ca^{2}%
}{\left( 1-\bar{t}\right) ^{p}}e^{\delta \left( t\right) }. 
\]

Now we turn to the lower bounds estimates. Since they are similar to those
given in detail in section 3, we obtain that the estimate $(4.24)$ for
potential operator function $\tilde{V}\left( x,t\right) $ when 
\[
\frac{\theta n}{2}-1<\frac{p}{2\left( 2-p\right) }\text{, i.e., }p>\frac{%
2\left( \theta n-2\right) }{\theta n-1}. 
\]

Finally, we get 
\[
e^{\frac{C}{\left( 1-\bar{t}\right) ^{p}}\log \phi \left( t\right) }\leq
e^{C\gamma R^{p}}\leq e^{\frac{C}{\left( 1-\bar{t}\right) ^{p/2}}\vartheta
\left( t\right) },\text{ }\vartheta \left( t\right) =C\left( 1-\bar{t}%
\right) ^{\frac{-p^{2}}{2\left( 2-p\right) }} 
\]%
for $p<\frac{p}{2}+\frac{p^{2}}{2\left( 2-p\right) }$, i.e. $p>1$ that
assumed in conditon of Theorem 4, i.e. we obtain the assertion of Theorem 4.

\textbf{Remark 5.1. }Let us consider the case $\theta =4/n$ in Theorem $4$,
i.e. $p>4/3$. Then from Theorem 4 we obtain the following result

\textbf{Result 5.1}. Assume the conditions of Theorem 4 are satisfied for $%
p>4/3$. Then $u\left( x,t\right) \equiv 0.$

\textbf{6. Unique continuation properties for the system of Schr\"{o}dinger
equation }

Consider the Cauchy problem for the finite or infinite system of Schr\"{o}%
dinger equation%
\begin{equation}
\frac{\partial u_{m}}{\partial t}=i\left[ \Delta
u_{m}+\sum\limits_{j=1}^{N}a_{mj}u_{j}+\sum\limits_{j=1}^{N}b_{mj}u_{j}%
\right] ,\text{ }x\in R^{n},\text{ }t\in \left( 0,T\right) ,  \tag{6.1}
\end{equation}%
where $u=\left( u_{1},u_{2},...,u_{N}\right) ,$ $u_{j}=u_{j}\left(
x,t\right) ,$ $a_{mj}$ are complex numbers and $b_{mj}=b_{mj}\left(
x,t\right) $ are complex valued functions$.$ Let $l_{2}=l_{2}\left( N\right) 
$ and $l_{2}^{s}=l_{2}^{s}\left( N\right) $ (see $\left[ \text{23, \S\ 1.18}%
\right] $). Let $A$ be the operator in $l_{2}\left( N\right) $ defined by%
\[
\text{ }D\left( A\right) =\left\{ u=\left\{ u_{j}\right\} \in l_{2},\text{ }%
\left\Vert u\right\Vert _{l_{2}^{s}\left( N\right) }=\left(
\sum\limits_{j=1}^{N}2^{sj}\left\vert u_{j}\right\vert ^{2}\right) ^{\frac{1%
}{2}}<\infty \right\} , 
\]

\[
A=\left[ a_{mj}\right] \text{, }a_{mj}=g_{m}2^{sj},\text{ }s>0,\text{ }%
m,j=1,2,...,N,\text{ }N\in \mathbb{N} 
\]

\bigskip and 
\[
\text{ }D\left( V\left( x,t\right) \right) =\left\{ u=\left\{ u_{j}\right\}
\in l_{2}^{s}\right\} , 
\]

\[
V\left( x,t\right) =\left[ b_{mj}\left( x,t\right) \right] \text{, }%
b_{mj}\left( x,t\right) =g_{m}\left( x,t\right) 2^{sj},\text{ }%
m,j=1,2,...,N. 
\]

Let 
\[
X_{2}=L^{2}\left( R^{n};l_{2}\right) ,Y^{s,2}=H^{s,2}\left(
R^{n};l_{2}\right) . 
\]

\ From Theorem 1 we obtain the following result

\textbf{Theorem 6.1. }Assume there exist the constants $a_{0},$ $a_{1},$ $%
a_{2}>0$ such that for any $k\in \mathbb{Z}^{+}$ a solution $u\in C\left( %
\left[ 0,1\right] ;X_{2}\right) $ of $\left( 6.1\right) $ satisfy 
\[
\dint\limits_{R^{n}}\left\Vert u\left( x,0\right) \right\Vert
_{l_{2}}^{2}e^{2a_{0}\left\vert x\right\vert ^{p}}dx<\infty ,\text{ for }%
p\in \left( 1,2\right) , 
\]%
\[
\dint\limits_{R^{n}}\left\Vert u\left( x,1\right) \right\Vert
_{l_{2}}^{2}e^{2k\left\vert x\right\vert ^{p}}dx<a_{2}e^{2a_{1}k^{\frac{q}{%
q-p}}},\text{ }\frac{1}{p}+\frac{1}{q}=1. 
\]%
Moreover, there exists $M_{p}>0$ such that

\[
a_{0}a_{1}^{p-2}>M_{p}. 
\]

Then $u\left( x,t\right) \equiv 0.$

\ \textbf{Proof.} It is easy to see that $A$ is a symmetric operator in $%
l_{2}$ and other conditions of Theorem 1 are satisfied. Hence, from Teorem 1
we obtain the conculision.

\begin{center}
\textbf{7. Unique continuation properties for nonlinear anisotropic Schr\"{o}%
dinger equation }

\ \ \ \ \ \ \ \ \ \ \ \ \ \ \ \ \ \ \ \ \ \ \ \ \ \ \ \ \ \ \ \ \ \ \ \ \ \
\ \ \ \ \ \ \ \ \ \ 
\end{center}

The regularity property of BVP for elliptic equations\ were studied e.g. in $%
\left[ \text{1, 2}\right] $. Let $\Omega =R^{n}\times G$, $G\subset R^{d},$ $%
d\geq 2$ is a bounded domain with $\left( d-1\right) $-dimensional boundary $%
\partial G$. Let us consider the following problem

\begin{equation}
i\partial _{t}u+\Delta _{x}u+\sum\limits_{\left\vert \alpha \right\vert \leq
2m}a_{\alpha }\left( y\right) D_{y}^{\alpha }u\left( x,y,t\right) +F\left( u,%
\bar{u}\right) u=0,\text{ }  \tag{7.1}
\end{equation}%
\[
\text{ }x\in R^{n},\text{ }y\in G,\text{ }t\in \left[ 0,1\right] ,\text{ } 
\]

\begin{equation}
B_{j}u=\sum\limits_{\left\vert \beta \right\vert \leq m_{j}}\ b_{j\beta
}\left( y\right) D_{y}^{\beta }u\left( x,y,t\right) =0\text{, }x\in R^{n},%
\text{ }y\in \partial G,\text{ }j=1,2,...,m,  \tag{7.2}
\end{equation}%
where $a_{\alpha },$ $b_{j\beta }$ are the complex valued functions, $\alpha
=\left( \alpha _{1},\alpha _{2},...,\alpha _{n}\right) $, $\beta =\left(
\beta _{1},\beta _{2},...,\beta _{n}\right) ,$ $\mu _{i}<2m,$ $K=K\left(
x,y,t\right) $ and 
\[
D_{x}^{k}=\frac{\partial ^{k}}{\partial x^{k}},\text{ }D_{j}=-i\frac{%
\partial }{\partial y_{j}},\text{ }D_{y}=\left( D_{1,}...,D_{n}\right) ,%
\text{ }y=\left( y_{1},...,y_{n}\right) . 
\]

$\ $

\bigskip Let%
\[
\xi ^{\prime }=\left( \xi _{1},\xi _{2},...,\xi _{n-1}\right) \in R^{n-1},%
\text{ }\alpha ^{\prime }=\left( \alpha _{1},\alpha _{2},...,\alpha
_{n-1}\right) \in Z^{n},\text{ } 
\]%
\[
\text{ }A\left( y_{0},\xi ^{\prime },D_{y}\right) =\sum\limits_{\left\vert
\alpha ^{\prime }\right\vert +j\leq 2m}a_{\alpha ^{\prime }}\left(
y_{0}\right) \xi _{1}^{\alpha _{1}}\xi _{2}^{\alpha _{2}}...\xi
_{n-1}^{\alpha _{n-1}}D_{y}^{j}\text{ for }y_{0}\in \bar{G} 
\]%
\[
B_{j}\left( y_{0},\xi ^{\prime },D_{y}\right) =\sum\limits_{\left\vert \beta
^{\prime }\right\vert +j\leq m_{j}}b_{j\beta ^{\prime }}\left( y_{0}\right)
\xi _{1}^{\beta _{1}}\xi _{2}^{\beta _{2}}...\xi _{n-1}^{\beta
_{n-1}}D_{y}^{j}\text{ for }y_{0}\in \partial G 
\]

\textbf{Theorem 7.1}. Let the following conditions be satisfied:

\bigskip (1) $G\in C^{2}$, $a_{\alpha }\in C\left( \bar{G}\right) $ for each 
$\left\vert \alpha \right\vert =2m$ and $a_{\alpha }\in L_{\infty }\left(
G\right) $ for each $\left\vert \alpha \right\vert <2m$;

(2) $b_{j\beta }\in C^{2m-m_{j}}\left( \partial G\right) $ for each $j$, $%
\beta $ and $\ m_{j}<2m$, $\sum\limits_{j=1}^{m}b_{j\beta }\left( y^{\prime
}\right) \sigma _{j}\neq 0,$ for $\left\vert \beta \right\vert =m_{j},$ $%
y^{^{\shortmid }}\in \partial G,$ where $\sigma =\left( \sigma _{1},\sigma
_{2},...,\sigma _{n}\right) \in R^{n}$ is a normal to $\partial G$ $;$

(3) for $y\in \bar{G}$, $\xi \in R^{n}$, $\lambda \in S\left( \varphi
_{0}\right) $ for $0\leq \varphi _{0}<\pi $, $\left\vert \xi \right\vert
+\left\vert \lambda \right\vert \neq 0$ let $\lambda +$ $\sum\limits_{\left%
\vert \alpha \right\vert =2m}a_{\alpha }\left( y\right) \xi ^{\alpha }\neq 0$%
;

(4) for each $y_{0}\in \partial G$ local BVP in local coordinates
corresponding to $y_{0}$:%
\[
\lambda +A\left( y_{0},\xi ^{\prime },D_{y}\right) \vartheta \left( y\right)
=0, 
\]

\[
B_{j}\left( y_{0},\xi ^{\prime },D_{y}\right) \vartheta \left( 0\right)
=h_{j}\text{, }j=1,2,...,m 
\]%
has a unique solution $\vartheta \in C_{0}\left( \mathbb{R}_{+}\right) $ for
all $h=\left( h_{1},h_{2},...,h_{n}\right) \in \mathbb{C}^{n}$ and for $\xi
^{\prime }\in R^{n-1};$

(5) there exist positive constants $b_{0}$ and $\theta $ such that a
solution $u\in C\left( \left[ -1,1\right] ;X_{2}\right) $ of $\left(
7.1\right) -\left( 7.2\right) $ satisfied 
\[
\left\Vert F\left( u,\bar{u}\right) \right\Vert _{L^{2}\left( G\right) }\leq
b_{0}\left\Vert u\right\Vert _{L^{2}\left( G\right) }^{\theta }\text{ for }%
\left\Vert u\right\Vert _{L^{2}\left( G\right) }>1; 
\]

(6) Suppose 
\[
\left\Vert u\left( .,t\right) \right\Vert _{L^{2}\left( R^{n}\times G\right)
}=\left\Vert u\left( .,0\right) \right\Vert _{L^{2}\left( R^{n}\times
G\right) }=\left\Vert u_{0}\right\Vert _{L^{2}\left( R^{n}\times G\right)
}=a 
\]%
for $t\in \left[ -1,1\right] $ and that $\left( 2.15\right) $ holds with $%
Q\left( .\right) $ satisfies $\left( 2.16\right) $ for $H=L^{2}\left(
G\right) .$

\ If $p>p\left( \theta \right) =\frac{2\left( \theta n-2\right) }{\left(
\theta n-1\right) },$ then $a\equiv 0.$

\ \textbf{Proof. }Let us consider operators $A$ and $V\left( x,t\right) $ in 
$H=L^{2}\left( G\right) $ that are defined by the equalities 
\[
D\left( A\right) =\left\{ u\in W^{2m,2}\left( G\right) \text{, }B_{j}u=0,%
\text{ }j=1,2,...,m\text{ }\right\} ,\ Au=\sum\limits_{\left\vert \alpha
\right\vert \leq 2m}a_{\alpha }\left( y\right) D_{y}^{\alpha }u\left(
y\right) , 
\]

Then the problem $\left( 7.1\right) -\left( 7.2\right) $ can be rewritten as
the problem $\left( 2.12\right) $, where $u\left( x\right) =u\left(
x,.\right) ,$ $f\left( x\right) =f\left( x,.\right) $,\ $x\in R^{n}$ are the
functions with values in\ $H=L^{2}\left( G\right) $. By virtue of $\left[ 
\text{1}\right] $ operator $A+\mu $ is positive in $L^{2}\left( G\right) $
for sufficiently large $\mu >0$. Moreover, in view of (1)-(6) all conditons
of Theorem 3 are hold. Then Theorem 3 implies the assertion.

\begin{center}
\textbf{8.} \textbf{The Wentzell-Robin type mixed problem for Schr\"{o}%
dinger equations}
\end{center}

Consider the problem $\left( 1.5\right) -\left( 1.6\right) $. \ Let 
\[
\sigma =R^{n}\times \left( 0,1\right) \text{, }X_{2}=L^{2}\left( R^{n}\times
\left( 0,1\right) \right) ,\text{ }Y^{2,k}=H^{2,k}\left( R^{n}\times \left(
0,1\right) \right) . 
\]

Suppose $\nu =\left( \nu _{1},\nu _{2},...,\nu _{n}\right) $ are nonnegative
real numbers. In this section, from Theorem 1 we obtain the following result:

\bigskip \textbf{Theorem 8.1. } Suppose the the following conditions are
satisfied:

(1)\ $a$ is positive, $b$ is a real-valued functions on $\left( 0,1\right) $%
. Moreover$,$ $a\left( .\right) \in C\left( 0,1\right) $ and%
\[
\exp \left( -\dint\limits_{\frac{1}{2}}^{x}b\left( t\right) a^{-1}\left(
t\right) dt\right) \in L_{1}\left( 0,1\right) ;
\]

(2) $u_{1},$ $u_{2}\in C\left( \left[ 0,1\right] ;Y^{2,k}\right) $ are
strong solutions of $(1.2)$ with $k>\frac{n}{2};$

(3) $F:\mathbb{C}\times \mathbb{C}\rightarrow \mathbb{C},$ $F\in C^{k}$, $%
F\left( 0\right) =\partial _{u}F\left( 0\right) =\partial _{\bar{u}}F\left(
0\right) =0;$

(4) there exist positive constants $\alpha $ and $\beta $ such that 
\[
e^{\frac{\left\vert \alpha x\right\vert ^{p}}{p}}\left( u_{1}\left(
.,0\right) -u_{2}\left( .,0\right) \right) \in X_{2},\text{ }e^{^{\frac{%
\left\vert \beta x\right\vert ^{q}}{q}}}\left( u_{1}\left( .,0\right)
-u_{2}\left( .,0\right) \right) \in X_{2}, 
\]

with 
\[
p\in \left( 1,2\right) \text{, }\frac{1}{p}+\frac{1}{q}=1; 
\]

(5) there exists $N_{p}>0$ such that 
\[
\alpha \beta >N_{p}.\text{ } 
\]

Then $u_{1}\left( x,t\right) \equiv u_{2}\left( x,t\right) .$

\ \textbf{Proof.} Let $H=L^{2}\left( 0,1\right) $ and $A$ is a operator
defined by $\left( 1.4\right) .$ Then the problem $\left( 1.5\right) -\left(
1.6\right) $ can be rewritten as the problem $\left( 1.2\right) $. By virtue
of $\left[ \text{10, 11}\right] $ the operator $A$ generates analytic
semigroup in $L^{2}\left( 0,1\right) $. Hence, by virtue of (1)-(5) all
conditons of Theorem 1.2 are satisfied. Then Theorem 1.2 implies the
assertion.

\textbf{References}\ \ 

\begin{enumerate}
\item H. Amann, Linear and quasi-linear equations,1, Birkhauser, Basel
(1995).

\item A. Benedek, A. Calder\`{o}n, R. Panzone, Convolution operators on
Banach space valued functions, Proc. Nat. Acad. Sci. USA, 48, 356--365 (1962)

\item A. Bonami, B. Demange, A survey on uncertainty principles related to
quadratic forms, Collect. Math., V. Extra, (2006), 1--36.

\item A.-P. Calder\'{o}n, Commutators of singular integral operators, Proc.
Nat. Acad. Sci. U.S.A. 53(1965), 1092--1099.

\item R. Denk, T. Krainer, $R$-Boundedness, pseudodifferential operators,
and maximal regularity for some classes of partial differential operators,
Manuscripta Math. 124(3) (2007), 319-342.

\item L. Escauriaza, C. E. Kenig, G. Ponce, L. Vega, On uniqueness
properties of solutions of Schr\"{o}dinger Equations, Comm. PDE. 31, 12
(2006) 1811--1823.

\item L. Escauriaza, C. E. Kenig, G. Ponce, and L. Vega, Uncertainty
principle of Morgan type and Schr\"{o}dinger evolution, J. London Math. Soc.
81, (2011) 187--207.

\item N. Hayashi, K. Nakamitsu, and N. Tsutsumi, On solutions of the initial
value problem for the nonlinear Schr%
\"{}%
odinger equations, J. Funct. Anal., 71 (1987).

\item L. H\"{o}rmander, A uniqueness theorem of Beurling for Fourier
transform pairs, Ark. Mat. 29, 2 (1991) 237--240.

\item J. A. Goldstain, Semigroups of Linear Operators and Applications,
Oxford University Press, Oxfard (1985).

\item A. Favini, G. R. Goldstein, J. A. Goldstein and S. Romanelli,
Degenerate Second Order Differential Operators Generating Analytic
Semigroups in $L_{p}$ and $W^{1,p}$, Math. Nachr. 238 (2002), 78 --102.

\item V. Keyantuo, M. Warma, The wave equation with Wentzell--Robin boundary
conditions on Lp-spaces, J. Differential Equations 229 (2006) 680--697.

\item S. G. Krein, Linear Differential Equations in Banach space, American
Mathematical Society, Providence, (1971).

\item A. Lunardi, Analytic Semigroups and Optimal Regularity in Parabolic
Problems, Birkhauser (2003).

\item J-L. Lions, E. Magenes, Nonhomogenous Boundary Value Broblems, Mir,
Moscow (1971).

\item V. B. Shakhmurov, Nonlinear abstract boundary value problems in
vector-valued function spaces and applications, Nonlinear Anal-Theor., 67(3)
2006, 745-762.

\item V. B. Shakhmurov, Hardy's uncertainty principle and unique
continuation properties for abstract Schr\"{o}dinger equations, Arkhif.
(2016).

\item R. Shahmurov, On strong solutions of a Robin problem modeling heat
conduction in materials with corroded boundary, Nonlinear Anal., Real World
Appl., v.13, (1), 2011, 441-451.

\item R. Shahmurov, Solution of the Dirichlet and Neumann problems for a
modified Helmholtz equation in Besov spaces on an annuals, J. Differential
equations, v. 249(3), 2010, 526-550.

\item C. Segovia, J. L.Torrea, Vector-valued commutators and applications,
Indiana Univ. Math.J. 38(4) (1989), 959--971.

\item E. M. Stein, R. Shakarchi, Princeton, Lecture in Analysis II. Complex
Analysis, Princeton University Press (2003).

\item A. Sitaram, M. Sundari, S. Thangavelu, Uncertainty principles on
certain Lie groups, Proc. Indian Acad. Sci. Math. Sci. 105 (1995), 135-151.

\item B. Simon, M. Schechter, Unique Continuation for Schrodinger Operators
with unbounded potentials, J. Math. Anal. Appl., 77 (1980), 482-492.

\item H. Triebel, Interpolation theory, Function spaces, Differential
operators, North-Holland, Amsterdam (1978).

\item S. Yakubov and Ya. Yakubov, Differential-operator Equations. Ordinary
and Partial \ Differential Equations, Chapman and Hall /CRC, Boca Raton
(2000).

\item Y. Xiao and Z. Xin, On the vanishing viscosity limit for the 3D
Navier-Stokes equations with a slip boundary condition. Comm. Pure Appl.
Math. 60, 7 (2007), 1027--1055.
\end{enumerate}

\end{document}